\documentclass[12pt]{article}

\begin{document}
\def\i{\indent}
\def\ni{\noindent}
\def\g{\gamma}
\def\a{\alpha}
\def\l{\lambda}
\def\m{\mu}
\def\s{\sigma}
\def\t{\theta}
\def\o{\omega}
\def\O{\Omega}
\def\e{\epsilon}
\def\H{\cal H}
\def\P{\Phi}
\def\p{\phi}
\def\b{\beta}
\begin{center}
{\bf SZEG\"O'S THEOREM AND ITS PROBABILISTIC DESCENDANTS}
\end{center}
\begin{center}
{\bf N. H. BINGHAM}
\end{center}
\begin{center}
{\bf Abstract}
\end{center}
The theory of orthogonal polynomials on the unit circle (OPUC)
dates back to Szeg\"o's work of 1915-21, and has been given a
great impetus by the recent work of Simon, in particular his
two-volume book [Si4], [Si5], the survey paper (or summary of
the book) [Si3], and the book [Si9], whose title we allude to in ours.  Simon's
motivation comes from spectral theory and analysis.  Another major
area of application of OPUC comes from probability, statistics,
time series and prediction theory; see for instance the book by
Grenander and Szeg\"o [GrSz].  Coming to the subject from this
background, our aim here is to complement [Si3] by giving some
probabilistically motivated results.  We also advocate a new definition of long-range dependence. \\

\ni {\it AMS 2000 subject classifications.} Primary 60G10, secondary 60G25. \\

\ni {\it Key words and phrases.} Stationary process, prediction theory, orthogonal polynomials
on the unit circle, partial autocorrelation function, moving average, autoregressive, long-range dependence, Hardy
space, cepstrum. \\

\begin{center}
{CONTENTS}
\end{center}
\S 1.  Introduction \\
\S 2.  Verblunsky's theorem and partial autocorrelation \\
\S 3.  Weak conditions: Szeg\"o's theorem \\
\S 4.  Strong conditions: Baxter's theorem \\
\S 5.  Strong conditions: the strong Szeg\"o theorem \\
\S 6.  Intermediate conditions \\
\i 6.1. Complete regularity \\
\i 6.2. Positive angle: the Helson-Szeg\"o and Helson-Sarason conditions \\
\i 6.3. Pure minimiality \\
\i 6.4. Rigidity; $(LM)$, $(CND)$, $(IPF)$ \\
\S 7.  Remarks \\
Acknowledgements \\
References \\

\ni {\bf \S 1.  Introduction}\\
\i The subject of orthogonal polynomials on the real line (OPRL),
at least some of which forms part of the standard undergraduate
curriculum, has its roots in the mathematics of the 19th century.
The name of Gabor Szeg\"o (1895-1985) is probably best remembered
nowadays for two things: co-authorship of 'P\'olya and Szeg\"o'
[PoSz] and authorship of 'Szeg\"o' [Sz4], his book of 1938, still
the standard work on OPRL.  Perhaps the key result in OPRL
concerns the
central role of the three-term recurrence relation ([Sz4], III.3.2: 'Favard's theorem'). \\
\i Much less well known is the subject of orthogonal polynomials on the unit circle (OPUC),
which dates from two papers of
Szeg\"o in 1920-21 ([Sz2], [Sz3]), and to which the last chapter of [Sz4] is devoted.  Again,
the key is the appropriate three-term recurrence relation, the {\it Szeg\"o recursion} or
{\it Durbin-Levinson algorithm} (\S 2).  This involves a sequence of coefficients (not
{\it two} sequences, as with OPRL), the {\it Verblunsky coefficients} $\a = ({\a}_n)$ (\S 2),
named (there are several other names in use) and systematically exploited in the magisterial
two-volume book on OPUC ([Si4], [Si5]) by Barry Simon.  See also his survey paper [Si3], written from the point of view of analysis and
spectral theory, the survey [GoTo], and his recent book [Si9].\\
\i Complementary to this is our own viewpoint, which comes from
probability and statistics, specifically {\it time series} (as
does the excellent survey of 1986 by Bloomfield [Bl3]). Here we
have a stochastic process (random phenomenon unfolding with time)
$X =(X_n)$ with $n$ integer (time discrete, as here, corresponds
to compactness of the unit circle by Fourier duality, whence the
relevance of OPUC; continuous time is also
important, and corresponds to OPRL).  \\
\i We make a simplifying assumption, and restrict attention to the {\it stationary}
case.  The situation is then invariant under the {\it shift} $n
\mapsto n+1$, which makes available the powerful mathematical
machinery of Beurling's work on invariant subspaces ([Beu]; [Nik1]).  While this is very
convenient mathematically, it is important to realize that this is both a strong restriction
and one unlikely to be satisfied exactly in practice.  One of the great contributions of the
statistician and econometrician Sir Clive Granger (1934-2009) was to demonstrate that
statistical/econometric methods appropriate for stationary situations can, when applied
indiscriminately to non-stationary situations, lead to misleading conclusions (via the
well-known statistical problem of {\it spurious regression}).  This has profound implications
for macroeconomic policy.  Governments depend on statisticians and econometricians for advice
on interpretation of macroeconomic data.  When this advice is misleading and mistaken policy
decisions are implemented, avoidable economic losses (in terms of GDP) may result which are
large-scale and permanent (cf. Japan's 'lost decade' in the 1990s, or lost two decades, and the global problems of
2007-8 on).\\
\i The mathematical machinery needed for OPUC is {\it function theory on the (unit) disc},
specifically the theory of Hardy spaces and Beurling's theorem (factorization into inner and
outer functions and Blaschke products).  We shall make free use of this, referring for what we need to standard works
(we recommend [Du], [Ho], [Gar], [Koo1], [Nik1], [Nik2]), but giving detailed references.  The theory on the disc (whose boundary the circle is compact)
corresponds analytically to the theory on the upper half-plane, whose boundary the real line is non-compact (for which see e.g. [DymMcK]).
Probabilistically, we work on the disc in discrete time and the half-plane in continuous time.  In each case, what dominates is an
{\it integrability condition}.  In discrete time, this is
{\it Szeg\"o's condition} (Sz), or {\it non-determinism} (ND) -- integrability of the
logarithm $\log \ w$ of the spectral density $w$ (of $\mu$) (\S 3).  In continuous time, this is the {\it logarithmic integral},
which gives its name to Koosis' book [Koo2]. \\
\i In view of the above, the natural context in which to work is that of complex-valued stochastic processes, rather than real-valued ones, in discrete time.
We remind the reader that here the Cauchy-Schwarz inequality tells us that correlation coefficients lie in the unit disc, rather than the interval $[-1,1]$. \\
\i The time-series aspects here go back at least as far as the work of Wiener [Wi1] in 1932 on
generalized harmonic analysis, GHA (which, incidentally, contains a good historical account of
the origins of spectral methods, e.g. in the work of Sir Arthur Schuster in the 1890s on
heliophysics).  During World War II, the {\it linear filter} (linearity is intimately linked
with Gaussianity) was developed independently by Wiener in the USA [Wi2], motivated by problems
of automatic fire control for anti-aircraft artillery, and Kolmogorov in Russia (then USSR) [Kol].  This work was developed by the Ukrainian mathematician M. G. Krein
over the period 1945-1985 (see e.g. [Dym]), by Wiener in the 1950s ([Wi3], IG, including commentaries) and by I. A. Ibragimov (1968 on). \\
\i The subject of time series is of great practical importance (e.g. in econometrics), but suffered within statistics by being regarded as 'for experts only'.  This changed with the 1970 book by Box and Jenkins (see [BoxJeRe]), which
popularized the subject by presenting a simplified account (including an easy-to-follow
model-fitting and model-checking recipe), based on ARMA models (AR for autoregressive, MA for
moving average).  The ARMA approach is still important; see e.g. Brockwell and Davis [BroDav] for a modern textbook account.  The realization that the Verblunsky coefficients $\a$ of OPUC are
actually the {\it partial autocorrelation function} (PACF) of time series opened the way for
the systematic exploitation of OPUC within time series by a number of authors.  These include
Inoue, in a series of papers from 2000 on (see especially [In3] of 2008), and Inoue and Kasahara from 2004 on (see especially [InKa2] of 2006). \\
\i Simon's work ([Si3], [Si4], [Si5]) focusses largely on four conditions, two weak (and
comparable) and two strong (and non-comparable).  Our aim here is to complement the expository
account in [Si3] by adding the time-series viewpoint.  This necessitates adding (at least) five
new conditions.  Four of these (comparable) we regard as intermediate, the fifth as strong.  In our view, one needs {\it three} levels of strength here, not two.  One is reminded of the Goldilocks principle (from the English children's story: not too hot/hard/high/..., not too cold/soft/low/..., but just right). \\
\i We begin in \S 2 by presenting the basics (Verblunsky's
theorem, PACF).  We turn in \S 3 to weak conditions (Szeg\"o's
condition (Sz), or (ND); Szeg\"o's theorem; $\a \in {\ell}_2$; $\sigma > 0$). In
\S 4 we look at our first strong condition, Baxter's condition
(B), and Baxter's theorem ($\a \in {\ell}_1$).  The satisfaction or otherwise of Baxter's condition (B) marks the transition between short- and long-range dependence.  The second strong
condition, the strong Szeg\"o condition (sSz), follows in \S 5
(strong Szeg\"o limit theorem, Ibragimov's theorem,
Golinskii-Ibragimov theorem, Borodin-Okounkov formula; $\a \in
H^{1/2}$), together with a weakening of (sSz), absolute
regularity.  We turn in \S 6 to intermediate conditions: in
decreasing order of strength, (i) complete regularity; (ii) positive angle (Helson-Szeg\"o, Helson-Sarason
and Sarason theorems); (iii) (pure) minimality (Kolmogorov); (iv) rigidity (Sarason),
Levinson-McKean condition (LM), complete non-determinism (CND), intersection of past and future
(IPF); see [KaBi] for details.  We close in \S 7 with some remarks. \\
\i The (weak) Szeg\"o limit theorem dates from 1915 [Sz1], the
strong Szeg\"o limit theorem from 1952 [Sz5].  Simon ([Si4], 11)
rightly says how remarkable it is for one person to have made
major contributions to the same area 37 years apart.  We note that
Szeg\"o's remarkable longevity here is actually exceeded (over the
40 years 1945-1985) by that of the late, great Mark Grigorievich
Krein (1907-1989). \\
\i What follows is a survey of this area, which contains (at least) eight different layers, of
increasing (or decreasing) generality.  This is an increase on Simon's (basic minimum of)
four.  We hope that no one will be deterred by this increase in
dimensionality, and so in apparent complexity.  Our aim is the precise opposite: to open up this fascinating area to a broader mathematical public, including the time-series, probabilistic and statistical communities.  For this, one needs to open up the `grey zone' between the strong and weak conditions, and examine the third category, of intermediate conditions .  We focus on these three levels of generality.  This largely reduces the effective dimensionality to three, which we feel simplifies matters.  Mathematics should be made as simple as possible, but not simpler (to adapt Einstein's immortal dictum about physics).\\
\i We close by quoting Barry Simon ([Si8], 85): "It's true that until Euclidean Quantum Field
Theory changed my tune, I tended to think of probabilists as a priesthood who translated
perfectly simple functional analytic ideas into a strange language that merely confused the
uninitiated."  He continues: in his 1974 book on Euclidean Quantum Field Theory, "the dedication
says: "To Ed Nelson who taught me how unnatural it is to view probability theory as unnatural" ". \\

\ni {\bf \S 2.  Verblunsky's theorem and partial autocorrelation}. \\
\i Let $X = (X_n: n \in Z)$ be a discrete-time, zero-mean,
(wide-sense) stationary stochastic process, with autocovariance
function $\g = ({\g}_n)$,
$$
{\g}_n = E[X_n \overline{X_0}]
$$
(the variance is constant by stationarity, so we may take it as 1, and then the autocovariance reduces to the autocorrelation). \\
\i Let $\cal H$ be the Hilbert space spanned by $X = (X_n)$ in the $L_2$-space of the
underlying probability space, with inner product $(X,Y) := E[X \overline{Y}]$ and norm
$\|X\| := [E(|X|^2)]^{1/2}$.  Write $T$ for the unit circle, the boundary of the unit disc $D$, parametrised by $z = e^{i \t}$; unspecified integrals are over $T$. \\

\ni {\bf Theorem 1 (Kolmogorov Isomorphism Theorem)}.  There is a process $Y$ on $T$ with orthogonal increments and a probability measure $\m$ on $T$ with \\
\ni (i)
$$
X_n = \int e^{int} dY(t);
$$
\ni (ii)
$$
E[dY(t)^2] = d\m(t).
$$
(iii) The autocorrelation function $\g$ then has the {\it spectral representation}
$$
{\g}_n = \int e^{-in\t} d \m(\t).
$$
(iv) One has the {\it Kolmogorov isomorphism} between $\cal H$ (the {\it time domain}) and $L_2(\mu)$ (the {\it frequency domain}) given by
$$
X_t \leftrightarrow e^{it.}, \eqno(KIT)
$$
for integer $t$ (as time is discrete).\\

\ni {\it Proof}. Parts (i), (ii) are the Cram\'er representation of 1942 ([Cra], [Do] X.4; Cram\'er and Leadbetter [CraLea] \S 7.5).  Part (iii), due originally to Herglotz in 1911,
follows from (i) and (ii)([Do] X.4, [BroDav] \S 4.3).  Part (iv) is due to Kolmogorov in 1941 [Kol].    All this rests on Stone's theorem of 1932, giving the spectral representation of groups of unitary transformations of linear operators
on Hilbert space; see [Do] 636-7 for a historical account and references (including work of Khintchine in 1934 in continuous time), [DunSch] X.5 for background on spectral theory. // \\

\i The reader will observe the link between the Kolmogorov Isomorphism Theorem and (ii), and its later counterpart from 1944, the It\^o Isomorphism Theorem and $(dB_t)^2 = dt$
in stochastic calculus. \\
\i To avoid trivialities, we suppose in what follows that $\m$ is {\it non-trivial} -- has infinite support.\\
 \i Since for integer $t$ the $e^{i t \theta}$ span polynomials in $e^{i \theta}$, prediction theory for stationary processes reduces to approximation by polynomials.  This is the classical approach to the main result of the subject, Szeg\"o's theorem (\S 2 below); see e.g. [GrSz], Ch. 3, [Ach],
Addenda, B.  We return to this in \S 7.7 below.\\
\i We write
$$
d\m(\t) = w(\t) d\t/2\pi + d{\m}_s(\t),
$$
so $w$ is the {\it spectral density} (w.r.t. normalized Lebesgue measure) and ${\m}_s$ is the
{\it singular part} of $\m$. \\
\i By stationarity,
$$
E[X_m \overline{X_n}] = {\g}_{|m-n|}.
$$
The {\it Toeplitz matrix} for $X$, or $\m$, or $\g$, is
$$
\Gamma := ({\g}_{ij}), \quad \hbox{where} \quad {\g}_{ij} :=
{\g}_{|i-j|}.
$$
It is positive definite. \\
\i For $n \in N$, write ${\H}_{[-n,-1]}$ for the subspace of $\H$
spanned by $\{X_{-n}, \ldots, X_{-1} \}$ (the finite past at time 0 of length $n$),
$P_{[-n,-1]}$
for projection onto ${\H}_{[-n,-1]}$ (thus $P_{[-n,-1]}X_0$ is the best linear predictor of
$X_0$
based on the finite past), $P_{[-n,-1]}^{\perp} := I - P_{[-n,-1]}$ for the orthogonal
projection (thus $P_{[-n,-1]}^{\perp}X_0 := X_0 - P_{[-n,-1]}X_0$ is the prediction error).  We
use a similar notation for prediction based on the infinite past.  Thus $\H_{(-\infty,-1]}$
is the closed linear span (cls) of $X_k$, $k \le -1$, $P_{(-\infty,-1]}$ is the corresponding
projection, and similarly for other time-intervals.  Write
$$
{\H}_n := {\H}_{(-\infty, n]}
$$
for the (subspace generated by) the past up to time $n$,
$$
{\H}_{-\infty} := \bigcap_{n = -\infty}^{\infty} {\H}_n
$$
for their intersection, the (subspace generated by) the {\it remote past}.  With $corr(Y,Z) := E[Y \overline{Z}]/\sqrt{E[|Y|^2].E[|Z|^2]}$ for $Y, Z$ zero-mean and not a.s. 0,  write also
$$
{\a}_n := corr(X_n - P_{[1,n-1]}X_n, X_0 - P_{[1,n-1]}X_0)
$$
for the correlation between the residuals at times $0$, $n$ resulting from (linear) regression
on the intermediate values $X_1, \ldots, X_{n-1}$.  The sequence
$$
\a = ({\a}_n)_{n=1}^{\infty}
$$
is called the {\it partial autocorrelation function} (PACF).  It is also called the sequence of {\it Verblunsky coefficients}, for reasons which will emerge below. \\

\ni {\bf Theorem 2 (Verblunsky's Theorem}.  There is a bijection between the sequences $\a = ({\a}_n)$ with each ${\a}_n \in D$ and the probability measures $\mu$ on $T$.\\

\i This result dates from Verblunsky in 1936 [V2], in connection with OPUC.  It was re-discovered long afterwards by Barndorff-Nielsen and Schou [BarN-S] in 1973 and Ramsey [Ram] in 1974, both in connection with parametrization of time-series models in statistics.  The Verblunsky bijection has the great advantage to statisticians of giving an
{\it unrestricted parametrization}: the only restrictions on the $\a_n$ are the obvious ones
resulting from their being correlations -- $|\a_n| \leq 1$, or as $\mu$ is non-trivial, $|\a_n| < 1$. By contrast,
$\g = (\g_n)$ gives a {\it restricted parametrization}, in that the possible values of $\g_n$ are
restricted by the inequalities of positive-definiteness (principal minors of the Toeplitz matrix
$\Gamma$ are positive).  This partly motivates the detailed study of the PACF in, e.g., [In1], [In2], [In3], [InKa1], [InKa2].  For general statistical background on partial autocorrelation, see e.g.
[KenSt], Ch. 27 (Vol. 2), \S 46.26-28 (Vol. 3).\\
\i As we mentioned in \S 1, the basic result for OPUC
corresponding to Favard's theorem for OPRL is {\it Szeg\"o's
recurrence} (or {\it recursion}): given a probability measure $\m$
on $T$, let
$\P_n$ be the monic orthogonal polynomials they generate (by
Gram-Schmidt orthogonalization). For every polynomial $Q_n$ of
degree $n$, write
$$
Q_n^{\ast}(z) := z^n \overline{Q_n(1/\bar z)}
$$
for the {\it reversed polynomial}.  Then the Szeg\"o recursion is
$$
\P_{n+1}(z) = z \P_n(z) - {\bar \a}_{n+1} \P_n^{\ast}(z),
$$
where the parameters $\a_n$ lie in $D$:
$$
|\a_n| < 1,
$$
and are the {\it Verblunsky coefficients} (also known variously as the Szeg\"o, Schur,
Geronimus and reflection coefficients; see [Si4], \S 1.1).  The double use of the name Verblunsky coefficients and the notation $\a = (\a_n)$ for the PACF and the coefficients is
justified: the two coincide.  Indeed, the Szeg\"o recursion is known in the time-series literature as the {\it Durbin-
Levinson algorithm}; see e.g. [BroDav], \S\S 3.4, 5.2.  The term Verblunsky coefficient is from Simon [Si4], to which we refer
repeatedly.  We stress that Simon writes $\a_n$ for our $\a_{n+1}$, and so has $n = 0,1,\ldots$ where we
have $n = 1,2,\ldots$.  Our notational convention is already established
in the time-series literature (see e.g. [BroDav], \S\S 3.4, 5.2), and is more convenient in our
context of the PACF, where $n = 1,2,\ldots$ has the direct interpretation as a time-lag between
past and future (cf. [Si4], (1.5.15), p. 56-57).  See [Si4], \S 1.5 and (for two proofs of
Verblunsky's theorem) \S 1.7, 3.1, and [McLZ] for a recent application of the unrestricted PACF
parametrization. \\
\i One may partially summarize the distributional aspects of Theorems 1 and 2 by the one-one correspondences
$$
\a \leftrightarrow \m \leftrightarrow \g.
$$

\ni {\it The Durbin-Levinson algorithm}\\
\i Write
$$
\hat X_{n+1} := \p_{n1} X_n + \ldots + \p_{nn} X_1
$$
for the best linear predictor of $X_{n+1}$ given $X_{n}, \ldots, X_{1}$,
$$
v_n := E[(X_{n+1} - \hat X_{n+1})^2] = E[(X_{n+1} - P_{[1,n]}X_{n+1})^2]
$$
for the mean-square error in the prediction of $X_{n+1}$ based on $X_1, \ldots,X_n$,
$$
{\bf \phi}_n := (\p_{n1}, \ldots, \p_{nn})^T \eqno(fpc)
$$
for the vector of {\it finite-predictor coefficients}.  The Durbin-Levinson algorithm ([Lev],
[Dur]; [BroDav] \S 5.2, [Pou] \S 7.2) gives the ${\bf \phi}_{n+1}$, $v_{n+1}$ {\it recursively}, in
terms of quantities known at time $n$, as follows:\\
(i) The first component of ${\bf \phi}_{n+1}$ is given by
$$
\p_{n+1,n+1} = [\g_{n+1} - \sum_{j=1}^n \p_{nj} \g_{n-j}]/v_n.
$$
The $\p_{nn}$ are the Verblunsky coefficients $\a_n$:
$$
\p_{nn} = \a_n.
$$
(ii) The remaining components are given by
$$
\pmatrix{\p_{n+1,1} \cr
         \vdots \cr
         \p_{n+1,n} \cr}
=
\pmatrix{\p_{n1} \cr
         \vdots \cr
         \p_{nn} \cr}
-\p_{n+1,n+1}\pmatrix{\p_{nn} \cr
         \vdots \cr
         \p_{n1} \cr}
         =
\pmatrix{\p_{n1} \cr
         \vdots \cr
         \p_{nn} \cr}
-\a_{n+1}\pmatrix{\p_{nn} \cr
         \vdots \cr
         \p_{n1} \cr}.
$$
(iii) The prediction errors are given recursively by
$$
v_0 = 1, \qquad v_{n+1} = v_n[1 - |\p_{n+1,n+1}|^2] = v_n[1 - |\a_{n+1}|^2].
$$
In particular, $v_n > 0$ and we have from (ii) that
$$
\p_{nj} - \p_{n+1,j} = \a_{n+1} \p_{n,n+1-j}. \eqno(DL).
$$
\ni Since by (iii)
$$
v_n = \prod_{j=1}^n [1 - |\a_n|^2],
$$
the $n$-step prediction error variance $v_n \to {\sigma}^2 > 0$ iff the infinite product
converges,
that is, $\a \in {\ell}_2$, an important condition that we will meet in \S 3 below in connection with Szeg\"o's condition.\\
{\it Note}. 1.  The Durbin-Levinson algorithm is related to the Yule-Walker equations of time-series
analysis (see e.g. [BroDav], \S 8.1), but avoids the need there for matrix inversion. \\
2.  The computational complexity of the Durbin-Levinson algorithm grows quadratically, rather
than cubically as one might expect; see e.g. Golub and van Loan [GolvL], \S 4.7.  Its good
numerical properties result from efficient use of the Toeplitz character of the matrix
$\Gamma$ (or equivalently, of Szeg\"o recursion). \\
3.  See [KatSeTe] for a recent approach to the Durbin-Levinson algorithm, and [Deg] for the multivariate case.\\

\ni {\it Stochastic versus non-stochastic} \\
\i This paper studies prediction theory for stationary stochastic processes.  As an extreme example (in which no prediction is possible),
take the `free' case, in which the $X_n$ are independent (and identically distributed).  Then $d \mu(\t) = d \t/2 \pi$, ${\g}_n = {\delta}_{n0}$, ${\a}_n \equiv 0$, ${\Phi}_n(z) = z^n$ ([Si4], Ex. 1.6.1). \\
\i In contrast to this is the situation where $X = (X_n)$ is {\it non-stochastic} -- deterministic, but (typically) chaotic.  This case often arises in non-linear time-series analysis and
dynamical systems; for a monograph treatment, see Kantz and Schreiber [KanSch]. \\

\i One natural way to classify results on OPUC is by the strength of
the conditions that they impose.  Simon's book discusses a range
of conditions, starting with a fairly weak one, Szeg\"o's
condition ([Si4] Ch. 2 and \S 3 below), and proceeding to two
principal stronger ones, Baxter's condition ([Si4] Ch. 5 and \S 4
below) and the strong Szeg\"o condition ([Si4] Ch. 6 and \S 5
below).  From a probabilistic viewpoint, equally important are a
range of intermediate conditions not discussed in Simon's book.  These we discuss in \S 6.
We close with some remarks in \S 7. \\

\ni {\bf \S 3. Weak conditions: Szeg\"o's theorem}.\\

\ni {\it Rakhmanov's Theorem} \\
\i One naturally expects that the influence of the distant past decays with increasing
lapse of time.  So one wants to know when
$$
\a_n \to 0 \qquad (n \to \infty).
$$
By {\it Rakhmanov's theorem} ([Rak]; [Si5] Ch 9, and Notes to \S 9.1, [MatNeTo]), this happens if
the density $w$ of the absolutely continuous component
$\mu_a$ is positive on a set of full measure:

$$
|\{ \theta: w(\theta) > 0 \}| = 1
$$
(using normalized Lebesgue measure -- or $2\pi$ using Lebesgue measure).\\

\ni{\it Non-determinism and the Wold decomposition}. \\
\i Write ${\s}^2$ for the one-step mean-square prediction error:
$$
{\s}^2 := E[(X_0 - P_{(-\infty, -1]} X_0)^2];
$$
by stationarity, this is the ${\s}^2 = \lim_{n \to \infty} v_n$ above.
Call $X$ {\it non-deterministic} $(ND)$ if $\s > 0$, {\it deterministic} if $\s = 0$.  (This
usage is suggested by the usual one of non-randomness being zero-variance, though here a
non-deterministic process may be random, but independent of time, so the stochastic
process reduces to a random variable.)  The {\it Wold decomposition} (von Neumann [vN] in 1929, Wold [Wo] in 1938;
see e.g. Doob [Do], XII.4, Hannan [Ha1], Ch. III) expresses a process $X$ as the sum of a non-deterministic process
$U$ and a deterministic process $V$:
$$
X_n = U_n + V_n;
$$
the process $U$ is a moving average,
$$
U_n = \sum_{j=-\infty}^n m_{n-j}{\xi}_j = \sum_{k=0}^{\infty} m_k {\xi}_{n-k},
$$
with the ${\xi}_j$ zero-mean and uncorrelated, with each other and with $V$; $E[{\xi}_n] = 0$, $var(\xi_n) = E[{\xi}_n^2] = {\s}^2$.  Thus when $\s = 0$ the ${\xi}_n$ are 0, $U$ is missing and the process is deterministic.  When $\s > 0$, the spectral measures of $U_n$, $V_n$ are $\m_{ac}$ and $\m_s$, the absolutely continuous and singular components of $\m$.  Think of ${\xi}_n$ as the `innovation' at time $n$ -- the new random input, a measure of the unpredictability of the present from the past.  This is only present when $\sigma > 0$; when $\sigma = 0$, the present is determined by the past -- even by the remote past.\\
\i The Wold decomposition arises in operator theory ([vN]; Sz.-Nagy and Foias in 1970 [SzNF], Rosenblum and Rovnyak in 1985 [RoRo], \S 1.3, [Nik2]), as a decomposition into the unitary and completely non-unitary (cnu) parts. \\

\ni {\it Szeg\"o's Theorem} \\

\ni {\bf Theorem 3 (Szeg\"o's Theorem)}. \\
(i) $\s > 0$ iff $\log w \in L_1$, that is,
$$
\int - \log w(\t) d \t > - \infty. \eqno(Sz)
$$
(ii) $\s > 0$ iff $\a \in {\ell}_2$. \\
(iii)
$$
{\s}^2 = {\prod}_1^{\infty} (1 - |{\a}_n|^2),
$$
so $\s > 0$ iff the product converges, i.e. iff
$$
{\sum} |{\a}_n|^2 < \infty: \qquad \a \in {\ell}_2;
$$
(iv) ${\s}^2$ is the geometric mean $G(\m)$ of $\m$:
$$
{\s}^2 = \exp \bigl(\frac{1}{2 \pi} \int \log w(\t) d\t \bigr) =: G(\m) > 0. \eqno(K)
$$

\ni {\it Proof}.  Parts (i), (ii) are due to Szeg\"o [Sz2], [Sz3] in 1920-21, with $\mu$ absolutely continuous, and to Verblunsky [V2] in 1936 for general $\mu$.
See [Si4] Ch. 2, [Si9] Ch. 2.  Parts (iii) and (iv) are due to Kolmogorov in 1941 [Kol].  Thus $(K)$ is called {\it Kolmogorov's formula}.
The alternative name for {\it Szeg\"o's condition} $(Sz)$ is the {\it non-determinism condition} $(ND)$, above. // \\

\i We now restrict attention to processes for which Szeg\"o's condition holds; indeed, we shall move below to stronger conditions.\\
\i The original motivation of Szeg\"o, and later Verblunsky, was {\it approximation theory}, specifically approximation by
polynomials.  The Kolmogorov Isomorphism Theorem allows us to pass between finite sections of the past to polynomials; denseness of polynomials allows prediction with zero error
(a `bad' situation -- determinism), which happens iff $(Sz)$ does not hold.  There is a detailed account of the (rather involved) history here in [Si4] \S 2.3.  Other classic contributions
include work of Krein in 1945, Levinson in 1947 [Lev] and Wiener in
1949 [Wi2].  See [BroDav] \S 5.8 (where un-normalized
Lebesgue measure is used, so there is an extra factor of $2 \pi$
on the right of $(K)$), [Roz] \S II.5 from the point of view of
time series, [Si4] for OPUC.\\

\ni {\it Pure non-determinism}, $(PND$ \\
\i When the remote past is trivial,
$$
{\H}_{-\infty} := \bigcap_{n = -\infty}^{\infty} {\H}_n = \{ 0 \}, \eqno(PND)
$$
there is no deterministic component in the Wold decomposition, and no singular component in the
spectral measure.  The process is then called {\it purely non-deterministic}.  Thus
$$
(PND) = (ND) + ({\m}_s = 0) = (Sz) + ({\m}_s = 0) = (\s > 0) + ({\m}_s = 0)  \eqno(PND)
$$
(usage differs here: the term `regular' is used for $(PND)$ in [IbRo], IV.1, but for
$(ND)$ in [Do], XII.2). \\

\ni {\it The Szeg\"o function and Hardy spaces} \\
\i Szeg\"o's theorem is the key result in the whole area, and to explore it further we need the Szeg\"o function ($h$, below).  For this, we need the language and viewpoint of the theory of Hardy spaces, and some of its standard results; several good textbook accounts are cited in \S 1.  For $0 < p < \infty$, the {\it Hardy space} $H_p$ is the class of analytic functions $f$ on $D$ for which
$$
{\sup}_{r < 1} \Bigl(\frac{1}{2 \pi} \int_0^{2 \pi} |f(r e^{i \t}|^p d \t \Bigr)^{1/p} < \infty. \eqno(H_p)
$$
As well as in time series and prediction, as here, Hardy spaces are crucial for martingale theory (see e.g. [Bin1] and the references there).  For an entertaining insight into Hardy spaces in
probability, see Diaconis [Dia]. \\
\i For non-deterministic processes, define the {\it Szeg\"o function} $h$ by
$$
h(z) := \exp \Bigl( \frac{1}{4 \pi} \int \Bigl(\frac{e^{i\t} + z}{e^{i\t} - z}
\Bigr) \log w(\t) d\t \Bigr) \qquad (z \in D), \eqno(OF)
$$
(note that in [In1-3], [InKa1,2], [Roz] II.5 an extra factor $\sqrt{2
\pi}$ is used on the right), or equivalently
$$
H(z) := h^2(z) = \exp \Bigl( \frac{1}{2 \pi} \int \Bigl(\frac{e^{i\t} + z}{e^{i\t} - z}
\Bigr) \log w(\t) d\t \Bigr) \qquad (z \in D).
$$
Because $\log w \in L_1$ by $(Sz)$, $H$ is an {\it outer function} for $H_1$ (whence the name
$(OF)$ above); see Duren [Du], \S 2.4.  By Beurling's canonical factorization theorem, \\
(i) $H \in H_1$, the Hardy space
of order 1 ([Du], \S 2.4), or as $H = h^2$, $h \in H_2$.\\
(ii) The radial limit
$$
H(e^{i\t}) := \lim_{r \uparrow 1} H(r e^{i\t})
$$
exists a.e., and
$$
|H(e^{i\t})| = |h(e^{i\t})|^2 = w(\t)
$$
(thus $h$ may be regarded as an `analytic square root' of $w$).
See also Hoffman [Ho], Ch. 3-5, Rudin [Ru], Ch. 17, Helson [He], Ch. 4.\\
\i Kolmogorov's formula now reads
$$
{\s}^2 = m_0^2 = h(0)^2 = G(\m) = \exp \bigl(\frac{1}{2 \pi} \int
\log w(\t) d\t \bigr).
 \eqno(K)
$$
When $\sigma > 0$, the Maclaurin coefficients $m = (m_n)$ of the Szeg\"o function $h(z)$ are the moving-average coefficients of the Wold decomposition (recall that the moving-average component does not appear when $\sigma = 0$); see Inoue [In3] and below.  When $\sigma > 0$, $m \in {\ell}_2$ is equivalent to convergence in mean square of the moving-average sum $\sum_{j=0}^{\infty} m_{n-j}{\xi}_j$ in the Wold decomposition.  This is standard theory for orthogonal expansions; see e.g. [Do], IV.4.  Note that a function being in $H_2$ and its Maclaurin coefficients being in ${\ell}_2$ are equivalent by general Hardy-space theory; see e.g. [Ru], 17.10 (see also Th. 17.17 for factorization), [Du] \S 1.4, 2.4, [Z2], VII.7.\\
\i Simon [Si4], \S 2.8 -- `Lots of equivalences' -- gives
Szeg\"o's theorem in two parts.  One ([Si4] Th. 2.7.14) gives
twelve equivalences, the other ([Si4], Th. 2.7.15) gives fifteen;
the selection of material is motivated by spectral theory [Si5].
Theorem 3 above extends these lists of equivalences, and treats
the material from the point of view of probability theory.  (It
does {\it not}, however, give a condition on the autocorrelation
$\g = ({\g}_n)$ equivalent to $(Sz)$; this is one of the outstanding problems of the area.) \\
\i The contrast here with Verblunsky's theorem is striking.  In
general, one has unrestricted parametrization: {\it all values}
$|{\a}_n|$ {\it are possible, for all} $n$.  But under Szeg\"o's
condition, one has $\a \in {\ell}_2$, and in particular ${\a}_n
\to 0$, as in Rakhmanov's theorem.  Thus non-deterministic processes
fill out only a tiny part of the
$\a$-parameter space $D^{\infty}$.  One may regard this as
showing that the remote past, trivial under $(Sz)$,
has a rich structure in general, as follows: \\

\ni {\it Szeg\"o's alternative} (or {\it dichotomy}).\\
\i  One either has
\begin{center}
$\log w \in L_1$ and ${\cal H}_{-\infty} \ne {\cal H}_{-n} \ne {\cal H}$
\end{center}
or
\begin{center}
$\log w \notin L_1$ and ${\cal H}_{-\infty} = {\cal H}_{-n} = {\cal H}$.
\end{center}
In the former case, $\a$ occupies a tiny part ${\ell}_2$ of $D^{\infty}$, and the remote past
${\cal H}_{-\infty}$ is identified with $L_2({\m}_s)$.  This is trivial iff ${\m}_s = 0$; cf.
$(PND)$.  In the second case, $\a$ occupies all of $D^{\infty}$, and the
remote past is the whole space. \\
\i Szeg\"o's dichotomy may be interpreted by analogy with physical systems.  Some systems
(typically, liquids and gases) are 'loose' -- left alone, they will thermalize, and tend to an
equilibrium in which the details of the past history are forgotten.  By contrast, some systems
(typically, solids) are 'tight': for example, in tempered steel, the thermal history is locked
in permanently by the tempering process.  Long memory is also important in economics and econometrics; for background here, see e.g. [Rob],
[TeKi]. \\
{\it Note}. 1. Our $h$ is the Szeg\"o function $D$ of Simon [Si4], (2.4.2), and $-1/h$ (see below) its negative reciprocal $-\Delta$ [Si4],
(2.2.92):
$$
h = D, \qquad -1/h = - \Delta
$$
(we use both notations to facilitate comparison between [In1-3], [InKa1,2], which use $h$, to
within the factor $\sqrt{2 \pi}$ mentioned above, and [Si4], our reference on OPUC, which
uses $D$). \\
2.  Both $h$ and $-1/h$ are analytic and non-vanishing in $D$.  See [Si4], Th. 2.2.14
(for $-1/h$, or $\Delta$), Th. 2.4.1 (for $h$, or $D$). \\
3.  That $(Sz)$ implies $h = D$ is in the unit ball of $H_2$ is in [Si4], Th. 2.4.1. \\
4.  See de Branges and Rovnyak [dBR] for general properties of such square-summable power
series.\\
5.  Our autocorrelation $\g$ is Simon's $c$ (he calls our
${\g}_n$, or his $c_n$, the {\it moments} of $\m$: [Si4],
(1.1.20)).  Our moving-average coefficients $m = (m_n)$ have no
counterpart in [Si4], and nor do the autoregressive coefficients
$r = (r_n)$ or minimality (see below for these).  We will
also need the Fourier coefficients of $\log w$, known for reasons
explained below as the {\it cepstrum}), which we write as $L =
(L_n)$ ('L for logarithm': Simon's ${\hat L}_n$ [Si4], (6.1.13)),
and a sequence $b = ({b}_n)$, the phase coefficients
(Fourier coefficients of ${\bar h}/h$). \\
6.  Lund et al. [LuZhKi] give several properties -- monotonicity, convexity etc. -- which one of
$m$, $\gamma$ has iff the other has. \\

\ni $MA(\infty)$ {\it and} $AR(\infty )$ \\
\i The power series expansion
$$
h(z) = \sum_{n=0}^{\infty} m_n z^n \qquad (z \in D)
$$
generates the $MA(\infty)$ coefficients $m = (m_n)$ in the Wold decomposition.  That of
$$
-1/h(z) = \sum_{n=0}^{\infty} r_n z^n \qquad (z \in D)
$$
generates the $AR(\infty)$ coefficients $r = (r_n)$ in the (infinite-order) {\it autoregression}
$$
\sum_{j = -\infty}^n r_{n-j} X_j + \xi_n = 0 \qquad (n \in Z). \eqno(AR)
$$
See [InKa2] \S 2, [In3] for background.\\
\i One may thus extend the above list of one-one correspondences, as follows:
$$
\hbox{Under $(Sz)$},\quad \a, \m, \g \leftrightarrow m = (m_n) \leftrightarrow h, -1/h
\leftrightarrow r = (r_n).
$$
\ni {\it Finite and infinite predictor coefficients}. \\
\i We met the $n$-vector ${\phi}_n$ of finite-predictor coefficients in $(fpc)$ of \S 1; we
can extend it to an infinite vector, still denoted ${\phi}_n$, by adding zeros.  The
corresponding vector $\phi := ({\phi}_1, {\phi}_2, \ldots)$ of {\it infinite-predictor
coefficients} gives the infinite predictor
$$
P_{(-\infty,-1]} X_0 = \sum_{j=1}^{\infty} {\phi}_j X_{-j}                    \eqno(ipc)
$$
([InKa2], (1.4)).  One would expect convergence of finite-predictor to infinite-predictor
coefficients; under Szeg\"o's condition, one has such convergence in ${\ell}_2$ iff
$(PND)$, i.e., ${\m}_s = 0$:
$$
{\phi}_n \to \phi \quad \hbox{in ${\ell}_2$} \Leftrightarrow (PND)
$$
(Pourahmadi [Pou], Th. 7.14).\\

\ni {\it The Szeg\"o limit theorem}.\\
\i With $G(\mu)$ as above, write $T_n$ (or $T_n(\g)$, or $T_n(\m)$)
for the $n \times n$ Toeplitz matrix $\Gamma^{(n)}$ with
elements
$$
\Gamma^{(n)}_{ij} := c_{j-i}
$$
obtained by truncation of the Toeplitz matrix $\Gamma$ (cf.
[BotSi2]).  {\it Szeg\"o's limit theorem} states that, under
$(Sz)$, its determinant satisfies
$$
\frac{1}{n} \log \det \ T_n \to G(\m) \qquad (n \to \infty)
$$
(note that $(Sz)$ is needed for the right to be defined).  A stronger statement --
{\it Szeg\"o's
strong limit theorem} -- holds; we defer this till \S 5.\\
\i The Szeg\"o limit theorem is used in the {\it Whittle estimator} of time-series analysis;
see e.g. Whittle [Wh], Hannan [Ha2].\\

\ni {\it Phase coefficients}.  \\
\i When the Szeg\"o condition $(Sz)$ holds, the Szeg\"o function
$h(z) = \sum_0^{\infty} m_n z^n$ is defined.  We can then define the {\it phase function}
${\bar h}/h$, so called because it has unit modulus and depends only on the phase or argument
of $h$ (Peller [Pel], \S 8.5).  Its Fourier coefficients $b_n$ are called the {\it phase
coefficients}.  They are given in terms of $m = (m_n)$ and $r = (r_n)$ by
$$
b_n := \sum_0^{\infty} m_k r_{n+k} \qquad (n = 0,1,2,\ldots). \eqno(b)
$$
The role of the phase coefficients is developed in [BiInKa].  They are important in connection
with rigidity (\S 6 below), and Hankel operators [Pel].\\

\ni{\it Rajchman measures}.  \\
\i In the Gaussian case, mixing in the
sense of ergodic theory holds iff
$$
\g_n \to 0 \qquad (n \to \infty)
$$
([CorFoSi], \S 14.2, Th. 2).  Since $(Sz)$ is $\g
\in {\ell}_2$, which implies $\g_n \to 0$, this is even weaker than $(Sz)$.  Measures for which this condition holds are called {\it
Rajchman measures} (they were studied by A. Rajchman in the
1920s).  Here the continuous singular part ${\m}_{cs}$ of $\mu$ is
decisive; for a characterization of Rajchman measures, see Lyons
([Ly1] -- [Ly3] and the appendix to [KahSa]).\\

\ni $ARMA(p,q)$.  \\
\i The Box-Jenkins $ARMA(p,q)$ methodology ([BoxJeRe],
[BroDav]: autoregressive of order $p$, moving average of order $q$ -- see \S 6.3 for $MA(q)$) applies
to stationary time series where the roots of the relevant polynomials lie in the unit disk (see e.g. [BroDav] \S 3.1).
The limiting case, of {\it unit roots}, involves non-stationarity, and so the statistical
dangers of spurious regression (\S 1); cf. Robinson [Rob], p.2.  We shall meet other instances
of unit-root phenomena later (\S 6.3). \\

\ni {\it Szeg\"o's theorem and the Gibbs Variational Principle} \\
\i We point out that Verblunsky [V2] proved the Gibbs Variational Principle, one of the cornerstones of nineteenth-century statistical mechanics, for the Szeg\"o integral:
$$
{\inf}_g [ \int e^g d \mu / \exp (\int g d \t /2 \pi )] = \exp [\int \log w(\t) d \t / 2 \pi].
$$
For details, see e.g. Simon [Si9] \S \S 2.2, 10.6, [Si10], Ch. 16, 17.  For background on the Gibbs Variational Principle, see e.g. Simon [Si1], III.4, Georgii [Geo], 15.4, Ellis [Ell], III.8. \\

\ni {\bf \S 4. Strong conditions: Baxter's theorem}\\

\i The next result ([Bax1], [Bax2], [Bax3]; [Si4], Ch. 5) gives the first of our strong conditions. \\

\ni {\bf Theorem 4 (Baxter's theorem)}.  The following are
equivalent: \\
(i) the Verblunsky coefficients (or PACF) are summable,
$$
\a \in {\ell}_1; \eqno(B)
$$
(ii) the autocorrelations are summable, $\g \in {\ell}_1$, and $\m$ is absolutely continuous
with continuous positive density:
$$
{\min}_{\t} w(\t) > 0.
$$

\i Of course, $({\g}_n)$ summable gives, as the ${\g}_n$ are the Fourier coefficients of $\m$,
that $\m$ is absolutely continuous with continuous density $w$; thus $w > 0$ iff $\inf w =
\min w > 0$.)  We extend this list of equivalences, and bring out its probabilistic significance, in Theorem 5 below on ${\ell}_1$
(this is substantially Theorem 4.1 of [In3]).  We call $\a \in {\ell}_1$ (or any of the other
equivalences in Theorem 4) {\it Baxter's condition} (whence $(B)$ above).  Since ${\ell}_1
\subset {\ell}_2$, Baxter's condition $(B)$ (`strong') implies Szeg\"o's condition $(Sz)$
(`weak').\\

\ni {\bf Theorem 5 (Inoue)}.  For a stationary process $X$, the following are equivalent: \\
(i) Baxter's condition $(B)$ holds: $\a \in {\ell}_1$.\\
(ii) $\g \in {\ell}_1$, ${\m}_s = 0$ and the spectral density $w$ is continuous and
positive.\\
(iii) $(PND)$ (that is, $(Sz)/(ND) + {\m}_s = 0$) holds, and the moving-average and
autoregressive coefficients are summable:
$$
m \in {\ell}_1, \qquad r \in {\ell}_1.
$$
(iv) $m \in {\ell}_1$, ${\m}_s = 0$ and the spectral density $w$ is continuous and
positive.\\
(v) $r \in {\ell}_1$, ${\m}_s = 0$ and the spectral density $w$ is continuous and
positive.\\

\ni {\it Proof}.\\
(i) $\Leftrightarrow$ (ii).  This is Baxter's theorem, as above. \\
(iii) $\Rightarrow$ (iv), (v).  By $(PND)$, $(Sz)$ holds, so the non-tangential limit
$$
h(e^{i\t}) = \lim_{r \uparrow 1} h(r e^{i\t}) = \sum_{n=0}^{\infty} m_n e^{in\t}
$$
exists a.e.  But as $m \in {\ell}_1$, $h(e^{i\t})$ is continuous, so this holds everywhere.
Since
$$
w(\t) = |h(e^{i\t})|^2 = |D(e^{i\t})|^2 = |\sum_{n=0}^{\infty} m_n e^{in\t}|^2,
$$
$w$ is continuous.  Letting $r \uparrow 1$ in
$$
h(z)(-1/h(z)) = (\sum_0^{\infty} m_n r^n e^{in\t})(\sum_0^{\infty} r_n r^n e^{in\t}) = -1
$$
gives similarly
$$
(\sum_0^{\infty} m_n e^{in\t})(\sum_0^{\infty} r_n e^{in\t}) = -1.
$$
So $h(e^{i\t})$ has no zeros, so neither
does $w$.  That is, (iv), (v) hold. \\
(iv) $\Rightarrow$ (iii).  As $w$ is positive and continuous, $w$ is bounded away from 0 and
$\infty$.  So $1/w$ is also.  So
$$
1/w(\t) = |1/h(e^{i\t})|^2 = |\Delta(e^{i\t})|^2
= |\sum_{n=0}^{\infty} r_n e^{in\t}|^2,
$$
where $\Delta = 1/D$.  (See [Si4], Th. 2.2.14, 2.7.15: the condition ${\lambda}_{\infty}(.) > 0$
there is $(Sz)$, so holds here.)  By Wiener's theorem, the reciprocal of a non-vanishing
absolutely convergent Fourier series is an absolutely convergent Fourier series (see e.g. [Ru],
Th. 18.21).  So from $m \in {\ell}_1$ we obtain $r \in {\ell}_1$, whence (iii)
(cf. [Berk], p.493).\\
(v) $\Rightarrow$ (iii).  This follows as above, by Wiener's theorem again. \\
(iv) $\Rightarrow$ (ii).  From the $MA(\infty)$ representation,
$$
{\g}_n = \sum_{k=0}^{\infty} m_{|n|+k} m_k \qquad (n \in Z) \eqno(conv)
$$
([InKa2], (2.21)).  So as ${\ell}_1$ is closed under convolution, $m \in {\ell}_1$ implies
$\g \in {\ell}_1$, indeed with
$$
\Vert \g \Vert_1 \le \Vert m \Vert_1^2,
$$
giving (ii). \\
(ii) $\Rightarrow$ (v).  We have
$$
{\phi}_j = c_0 r_j = \sigma r_j
$$
with ${\phi}_j$ the infinite-predictor coefficients ([InKa2], (3.1)).  Then $r \in
{\ell}_1$ follows by the Wiener-L\'evy theorem, as in Baxter [Ba3], 139-140. // \\

\ni {\it Note}. 1. Under Baxter's condition, both $|h|$ and $|1/h|$ (or $|D|$ and $|\Delta| = |1/D|$) are
continuous and positive on the unit circle.  As $h$, $1/h$ are analytic in the disk, and so
attain their maximum modulus on the circle by the maximum principle,
$$
\inf_D |h(.)| > 0, \qquad \inf_D |1/h(.)| > 0
$$
(and similarly for $D(.)$, $\Delta$); [Si4], (5.2.3), (5.2.4). \\
2.  The hard part of Baxter's theorem is (ii) $\Rightarrow$ (i), as Simon points out
([Si4], 314). \\
3.  Simon [Si4], Th. 5.2.2 gives {\it twelve} equivalences in his final form of Baxter's
theorem.  (He does not, however, deal explicitly with $m$ and $r$.) \\
4.  Simon also gives a more general form, in terms of {\it Beurling weights}, $\nu$.  The
relevant Banach algebras contain the Wiener algebra used above as the special case $\nu = 1$. \\
5.  The approach of [Si4], \S 5.1 is via truncated Toeplitz matrices and their inverses.  The
method derives, through Baxter's work, from the Wiener-Hopf technique.  This point of view is
developed at length in [BotSi1], [BotSi2].  Baxter's motivation was approximation to
infinite-past predictors by finite-past predictors. \\

\ni {\it Long-range dependence} \\
\i In various physical models, the property of {\it long-range dependence} (LRD) is important,
particularly in connection with phase transitions (see e.g. [Si1], Ch. II, [Gri1], Ch. 9,
[Gri2], Ch. 5), to which we return below.  This is a spatial
property, but applies also in time rather than space, when the term used is {\it long memory}.
A good survey of long-memory processes was given by Cox [Cox] in 1984, and a monograph
treatment by Beran [Ber] in 1994.  For more recent work, see [DouOpTa], [Rob], [Gao]
Ch. 6, [TeKi], [GiKoSu].\\
\i Baxter's theorem is relevant to the definition of LRD recently proposed independently
by Debowski [Deb] and Inoue [In3]: {\it long-range dependence, or long memory, is non-summability
of the PACF}:
\begin{center}
{$X$ has LRD iff $\a \notin {\ell}_1$.} \hfil $(DI)$ \break
\end{center}
\i While the broad concept of long memory, or LRD, has long been widely accepted, authors
differed over the
precise definition.  There were two leading candidates: \\
(i) LRD is non-summability of covariances, $\g \notin {\ell}_1$.\\
(ii) LRD is covariance decaying like a power: ${\g}_n \sim c/n^{1-2d}$ as $n \to \infty$, for
some parameter $d \in (0, 1/2)$ ($d$ for differencing -- see below) and constant
$c \in (0,\infty)$ (and so $\sum {\g}_n = \infty$).\\
{\it Note}. 1.  In place of (ii), one may require $w(\t) \sim C/{\t}^{2d}$ as $\t \downarrow 0$,
for some constant $C \in (0,\infty)$.  The constants here may be replaced by slowly varying
functions.  See e.g. [BinGT] \S 4.10 for relations between regular variation of Fourier series
and Fourier coefficients. \\
2.  One often encounters, instead of $d \in (0, 1/2)$, a parameter $H = d + \frac{1}{2} \in
(1/2, 1)$.  This $H$ is the {\it Hurst parameter}, named after the classic studies by the
hydrologist Hurst of water flows in the Nile; see [Ber], Ch. 2. \\
3.  For $d \in (0, \frac{1}{2})$, $\ell(.)$ slowly varying, the following class of prototypical
long-memory examples is considered in [InKa2], \S 2.3 (see also [In1], Th. 5.1):
$$
{\g}_n \sim \ell(n)^2 B(d,1-2d)/n^{1-2d},
$$
$$
m_n \sim \ell(n)/n^{1-d},
$$
$$
r_n \sim \frac{d \sin(\pi d)}{\pi}.\frac{1}{\ell(n)}.1/n^{1+d}.
$$
See the sources cited for inter-relationships between these. \\
4.  In [InKa2], Example 2.6, the class of $FARIMA(p,d,q)$ processes is considered (obtained from
an $ARIMA(p,q)$ process by fractional differencing of order $d$ -- see [Hos], [BroDav], [KokTa]).  For $d \in (0, 1/2)$ these
have long memory; for $d = 0$ they reduce to the familiar $ARMA(p,q)$ processes.  \\
\i Li ([Li], \S 3.4) has recently given a related but different definition of long
memory; we return to this in \S 5 below.\\

\ni {\bf 5.  Strong conditions: the strong Szeg\"o theorem}\\

\i The work of this section may be motivated by work from two areas of physics. \\

\ni {\it 1. The cepstrum}.\\
\i During the Cold War, the problem of determining the signature of the underground explosion
in a nuclear weapon test, and distinguuishing it from that of an earthquake, was very important, and was studied by
the American statistician J. W. Tukey and collaborators.  Write $L = (L_n)$, where the
$L_n$ are the Fourier coefficients of
$\log w$, the log spectral density:
$$
L_n := \int \log w(\t) e^{in\t} d\t/2 \pi.
$$
Thus $\exp(L_0)$ is the geometric mean $G(\m)$.  The sequence $L$ is called the {\it cepstrum},
$L_n$ the ceptstral coefficients (Simon's notation here is ${\hat L}_n$; [Si4], (2.1.14),
(6.1.11)); see e.g. [OpSc], Ch. 12.  The
terminology dates from work of Bogert, Healy and Tukey of 1963 on echo detection [BogHeTu]; see
McCullagh [McC], Brillinger [Bri] (the
term is chosen to suggest both echo and spectrum, by reversing the first half of the word
spectrum; it is accordingly pronounced with the c hard, like a k).\\

\ni {\it 2. The strong Szeg\"o limit theorem}.\\
\i This  (which gives the weak form on taking logarithms) states
(in its present form, due to Ibragimov) that
$$
\frac{\det \ T_n}{G(\m)^n} \to E(\m) := \exp \{ \sum_1^{\infty} k L_k^2) \} \qquad
(n \to \infty)
$$
(of course the sum here must converge; it turns out that this form is best-possible: the
result is valid whenever it makes sense ([Si4], 337).\\
\i The motivation was Onsager's work in the two-dimensional
Ising model, and in particular {\it Onsager's formula}, giving the existence of a critical
temparature $T_c$ and the decay of the magnetization as the temperature $T \uparrow T_c$;
see [BotSi2] \S 5.1, [Si1] II.6, [McCW].  The mechanism was a question by Onsager (c. 1950)
to his Yale colleague Kakutani, who asked Szeg\"o ([Si4], 331).\\

\i Write $H^{1/2}$ for the subspace of ${\ell}_2$ of sequences $a = (a_n)$ with
$$
{\Vert a \Vert}^2 := \sum_n (1 + |n|) |\a_n|^2 < \infty \eqno(H^{1/2})
$$
(the function of the `1' on the right is to give a norm; without it, $\Vert . \Vert$
vanishes on the constant functions).  This is a Sobolev space ([Si4], 329, 337; it is also a
Besov space, whence the alternative
notation $B_2^{1/2}$; see e.g. Peller [Pel], Appendix 2.6 and \S 7.13).  This is the space that
plays the role
here of ${\ell}_2$ in \S 2 and ${\ell}_1$ in \S 3.  Note first that, although ${\ell}_1$
and $H^{1/2}$ are close in that a sequence $(n^c)$ of powers belongs to both or
neither, neither contains the other (consider $a_n = 1/(n \log n)$, $a_n = 1/\sqrt{n}$ if
$n = 2^k$, 0 otherwise). \\

\ni {\bf Theorem 6 (Strong Szeg\"o Theorem)}. \\
(i)  If $(PND)$ holds (i.e. $(Sz) = (ND)$ holds and ${\m}_s = 0$), then
$$
E(\mu)
= \prod_{j=1}^{\infty} (1 - |{\a}_j|^2)^{-j}
= \exp \Bigl(\sum_{n=1}^{\infty} n L_n^2)
$$
(all three may be infinite), with the infinite product converging iff the {\it strong Szeg\"o condition}
$$
\a \in H^{1/2}, \eqno(sSz)
$$
holds.\\
(ii) $(sSz)$ holds iff
$$
L \in H^{1/2} \eqno(sSz')
$$
holds. \\
(iii) Under $(Sz)$, finiteness of any (all three) of the expressions in (i) forces ${\mu}_s = 0$. \\

\ni {\it Proof}.  Part (i) is due to Ibragimov ([Si4], Th. 6.1.1), and (ii) is immediate from this.  Part (iii) is due to
Golinski and Ibragimov ([Si4], Th. 6.1.2; cf. [Si2]). // \\

\i Part of Ibragimov's theorem was recently obtained independently by Li [Li], under the term
{\it reflectrum identity} (so called because it links the Verblunsky or reflection coefficients
with the cepstrum), based on information theory -- mutual information between past and
future.  Earlier, Li and Xie [LiXi] had shown the following: \\
(i) a process with given autocorrelations ${\gamma}_0, \ldots, {\gamma}_p$ with minimal
information between past and future must be an autoregressive model $AR(p)$ of order $p$; \\
(ii) a process with given cepstral coefficients $L_0, \ldots, L_p$ with minimal information
between past and future must be a {\it Bloomfield model} $BL(p)$ of order $p$ ([Bl1], [Bl2]),
that is, one with spectral density $w(\theta)
= \exp \{L_0 + 2 \sum_{k=1}^p L_k \cos \ k \theta \}$.\\
\i Another approach to the strong Szeg\"o limit theorem, due to Kac [Kac], uses the conditions
$$
\inf w(.) > 0, \qquad \g = ({\g}_n) \in {\ell}_1, \qquad \g \in H^{1/2}
$$
(recall that ${\ell}_1$ and $H^{1/2}$ are not comparable).  This proof, from 1954, is linked
to probability theory -- Spitzer's identity of 1956, and hence to fluctuation theory for random walks, for which see e.g. [Ch], Ch. 8.\\

\ni {\it The Borodin-Okounkov formula}.\\
\i This turns the strong Szeg\"o limit theorem above from analysis
to algebra by identifying the quotient on the left there as a
determinant which visibly tends to 1 as $n \to \infty$ [BorOk];
see [Si4] \S 6.2.  (It was published in 2000, having been previously
obtained by Geronimo and Case [GerCa] in 1979; see [Si4] 337, 344,
[Bot] for background here.)  In terms of operator theory and in Widom's notation [Bot], the
result is
$$
\frac{det \ T_n(a)}{G(a)^n}
 = \frac{det(I - Q_n H(b) H(\tilde c) Q_n)}{det (I - H(b)H(\tilde c))},
$$
for $a$ a sufficiently smooth function without zeros on the unit circle and with winding number
0. Then $a$ has a Wiener-Hopf factorization $a = a_- a_+$; $b := a_- a_+^{-1}$,
$c := a_-^{-1} a_+$; $H(b)$, $H(\tilde c)$ are the Hankel matrices
$H(b) = (b_{j+k+1})_{j,k = 0}^{\infty}$, $H(\tilde c) = (c_{-j-k-1})_{j,k = 0}^{\infty}$, and
$Q_n$ is the orthogonal projection of ${\ell}^2({1,2, \ldots})$ onto ${\ell}^2(\{n,n+1,\ldots \})$.
By Widom's formula,
$$
1/det(I - H(b)H(\tilde c)) = \exp \{ \sum_{k=1}^{\infty} k L_k^2 \} =: E(a)
$$
(see e.g. [Si4], Th. 6.2.13), and $Q_n H(b) H(\tilde c) Q_n \to 0$ in the trace norm, whence
$$
det \ T_n(a)/G(a)^n \to E(a),
$$
the strong Szeg\"o limit theorem.  See [Si4], Ch. 6, [Si6], [BasW], [BotW] (in [Si4] \S 6.2 the
result is given in OPUC terms; here $b$, $c$ are the phase function $\overline h/h$ and its
inverse).\\

\ni $(B + sSz)$.\\
\i We may have both of the strong conditions $(B)$ and $(sSz)$ (as happens in Kac's
method [Kac], for instance).  Matters then simplify, since the spectral density $w$ is now
continuous and positive.  So $w$ is bounded away from 0 and $\infty$, so $\log \ w$ is bounded.  Write
$$
{\omega}^2(\delta, h) :=
\sup_{|\t|\le \delta}\Bigl(\int |h(\lambda + \t) - h(\lambda)|^2 d \lambda \Bigr)^{1/2}
$$
for the $L_2$ modulus of continuity.  Applying [IbRo], IV.4, Lemma 7 to
$\log w$,
$$
L \in H^{1/2} \Leftrightarrow \sum_{k=1}^{\infty} {\omega}^2 (1/k, \log w) < \infty,
$$
and applying it to $w$,
$$
\g \in H^{1/2} \Leftrightarrow \sum_{k=1}^{\infty} {\omega}^2 (1/k, w) < \infty.
$$
Thus under $(B)$, $L \in H^{1/2}$ and $\g \in H^{1/2}$ become equivalent.  This last
condition is Li's proposed definition of long-range dependence:
$$
LRD \Leftrightarrow \g \notin H^{1/2} \eqno(Li)
$$
([Li], \S 3.4; compare the Debowski-Inoue definition $(DI)$ above, that LRD iff $\a \notin {\ell}_1$).\\
\i We are now in $W \cap H^{1/2}$, the intersection of $H^{1/2}$ with the Wiener algebra $W$
(of absolutely convergent Fourier series) relevant to Baxter's theorem as in \S 3.  As there,
we can take inverses, since the Szeg\"o function is non-zero on the circle (cf. [BotSi2], \S 5.1).
One can thus extend Theorem 2 to this situation, including the cepstral condition
$L \in H^{1/2}$ (Li [Li], Th. 1 part 3, showed that $L \in H^{1/2}$ and $\g \in H^{1/2}$ are
equivalent if $w$ is continuous and positive).\\

\ni $L_{\infty} + (sSz)$.\\
\i The bounded functions in $H^{1/2}$ form an algebra, the {\it Krein algebra} $K$, a Banach
algebra under convolution; see Krein [Kr], B\"ottcher and Silbermann [BotSi1] Ch. 10, [BotSi2]
Ch. 5, [Si4], 344, [BotKaSi].  The Krein algebra may be used as a partial substitute for the
Wiener algebra $W$used to treat Baxter's theorem in \S 3 ($W \cap H^{1/2}$ is also an algebra: [BotSi], \S 5.1).\\

\ni {\it 5.1}. $\phi$-{\it mixing}\\
\i Weak dependence may be studied by a hierarchy of mixing conditions; for background, see e.g.
Bradley [Bra1], [Bra2], [Bra3], Bloomfield [Bl3], Ibragimov and Linnik [IbLi], Ch. 17, Cornfeld et al. [CorFoSi]),
and in the Gaussian case Ibragimov and Rozanov [IbRo], Peller [Pel].  We need two sequences of mixing coefficients:
$$
\phi(n) := E \sup \{ |P(A | {\cal F}_{-\infty}^0 ) - P(A) | : A \in \hbox{${\cal F}_n^{\infty}$} \};
$$
$$
\rho(n) := \rho({\cal F}_{-\infty}^0, \hbox{${\cal F}_n^{\infty}$}),
$$
where
$$
\rho({\cal A}, \hbox{${\cal B}$}) := \sup \{ \Vert E(f | \hbox{${\cal B}$}) - E f \Vert_2/\Vert f \Vert_2: f \in L_2(\hbox{${\cal A}$}) \}.
$$
The process is called $\phi$-{\it mixing} if ${\phi}(n) \to 0$ as $n \to \infty$, $\rho$-{\it mixing} if $\rho(n) \to 0$.
(The reader is warned that some authors use other letters here -- e.g. [IbRo] uses $\beta$ for our $\phi$; we follow Bradley.)\\
\i We quote [Bra1] that $\phi$-mixing implies $\rho$-mixing.  We regard the first as a strong condition, so include it here, but the second and
its several weaker relatives as intermediate conditions, which we deal with in \S 6 below. \\
\i The spectral characterization for $\phi$-mixing is
$$
{\mu}_s = 0, \qquad w(\t) = |P(e^{i \t})|^2 w^{\ast}(\t),
$$
where $P$ is a polynomial with its roots on the unit circle and
the cepstrum $L^{\ast} = (L_n^{\ast})$ of $w^{\ast}$ satisfies the
strong Szeg\"o condition $(sSz)$ ([IbRo] IV.4, p. 129).  This is
weaker than $(sSz)$.  In the Gaussian case, $\phi$-mixing (also known as absolute regularity)
can also be characterized in operator-theoretic terms: $\phi(n)$ can be identified as
$\sqrt{tr(B_n)}$, where $B_n$ are compact operators with finite
trace, so $\phi$-mixing is $tr(B_n) \to 0$ ([IbRo], IV.2 Th. 4, IV.3 Th. 6).\\

\ni {\bf 6.  Intermediate conditions}\\

\ni We turn now to four intermediate conditions, in decreasing order of strength.\\

\ni {\it 6.1.  $\rho$-mixing} \\
\i The spectral characterization of $\rho$-mixing (also known as complete regularity) is
$$
\m_s = 0, \qquad w(\t) = |P(e^{i \t})|^2 w^{\ast}(\t),
$$
where $P$ is a polynomial with its roots on the unit circle and
$$
\log \ w^{\ast} = u + \tilde v,
$$
with $u$, $v$ real and continuous (Sarason [Sa2]; Helson and Sarason [HeSa]).   An alternative spectral characterization is
$$
\m_s = 0, \qquad w(\t) = |P(e^{i \t})|^2 w^{\ast}(\t),
$$
where $P$ is a polynomial with its roots on the unit circle and for all $\epsilon >0$,
$$
\log \ w^{\ast} = r_{\e} + u_{\e} + \tilde v_{\e},
$$
where $r_{\e}$ is continuous, $u_{\e}$, $v_{\e}$ are real and bounded, and
$\Vert u_{\e} \Vert + \Vert v_{\e} \Vert < \e$ ([IbRo], V.2 Th. 3; we note here that inserting
such a polynomial factor preserves complete regularity, merely changing $\rho$ -- [IbRo] V.1,
Th. 1).\\

\ni {\it 6.2. Positive angle: the Helson-Szeg\"o and Helson-Sarason conditions}.\\
\i We turn now to a weaker condition.  For subspaces $A$, $B$ of $\cal H$, the {\it angle}
between $A$ and $B$ is defined as
$$
{\cos}^{-1} \sup \{ |(a,b)|: a \in A, b \in B \}.
$$
Then $A$, $B$ are at a {\it positive angle} iff this supremum is $< 1$.  One says that the
process $X$ satisfies the {\it positive angle condition}, $(PA)$, if for some time lapse $k$ the past $cls(X_m: m < 0)$
and the future $cls(X_{k+m}: m \ge 0)$ are at a positive angle, i.e.
$\rho(0) = \ldots \rho(k-1) = 1, \rho(k)< 1$, which we write as $PA(k)$ (Helson and Szeg\"o
[HeSz], $k = 1$; Helson and Sarason [HeSa], $k > 1$).  The spectral characterization of this is
$$
\m_s = 0, \qquad w(\t) = |P(e^{i \t})|^2 w^{\ast}(\t),
$$
where $P$ is a polynomial of degree $k-1$ with its roots on the unit circle and
$$
\log \ w^{\ast} = u + \tilde v,
$$
where $u$, $v$ are real and bounded and $\Vert v \Vert < \pi/2$ ([IbRo] V.2, Th. 3, Th. 4).  (The
role of $\pi/2$ here stems from Zygmund's theorem of 1929, that if $u$ is bounded and
$\Vert u \Vert < \pi/2$, $\exp \{ \tilde u \} \in L_1$ ([Z1], [Z2] VII, (2.11), [Tor], V.3: cf.
[Pel] \S 3.2.)  Thus $\rho$-mixing implies $(PA)$ (i.e. $PA(k)$ for some $k$).\\
\i The case $PA(k)$ for $k > 1$ is a unit-root phenomenon (cf. the note at the end of \S 3).  We
may (with some loss of information) reduce to the case $PA(1)$ by sampling only at every $k$th
time point (cf. [Pel], \S \S 8.5, 12.8).  We shall do this for convenience in what follows. \\
\i It turns out that the Helson-Szeg\"o condition $(PA(1))$ coincides with {\it Muckenhoupt's
condition} $A_2$ in analysis:
$$
\sup_I \Bigl(\frac{1}{|I|} \int_I w(\t) d\t \Bigr)
\Bigl(\frac{1}{|I|} \int_I \frac{1}{w(\t)} d\t \Bigr) < \infty, \eqno(A_2)
$$
where $|.|$ is Lebesgue measure and the supremum is taken over all subintervals $I$ of the unit
circle $T$.  See e.g. Hunt, Muckenhoupt and Wheeden [HuMuWh].  With the above reduction of $PA$ to
$PA(1)$, we then have $\rho$-mixing implies $PA(1)$ (= $A_2)$. \\

\ni {\it 6.3.  Pure minimality} \\
\i Consider now the {\it interpolation problem}, of finding the best linear interpolation of a
missing value, $X_0$ say, from the others.  Write
$$
H_n' := cls \{X_m: m \ne n \}
$$
for the closed linear span of the values at times other than $n$.  Call $X$ {\it minimal} if
$$
X_n \notin H_n',
$$
{\it purely minimal} if
$$
\bigcap_n H_n' = \{ 0 \}.
$$
The spectral condition for minimality is (Kolmogorov in 1941, [Kol] \S 10)
$$
1/w \in L_1, \eqno(min)
$$
and for pure minimality, $\m_s = 0$ also (Makagon-Weron in 1976, [MakWe]; Sarason in 1978, [Sa1];
[Pou], Th. 8.10):
$$
1/w \in L_1, \qquad \m_s = 0. \eqno(purmin)
$$
Of course $(A_2)$ implies $1/w \in L_1$, so the Helson-Szeg\"o condition $(PA(1))$ (or
Muckenhoupt condition $(A_2)$) implies pure minimality.  (From $\log x < x-1$ for $x > 1$,
$(min)$ implies $(Sz)$: both restrict the
small values of $w \ge 0$, and in particular force $w > 0$ a.e.)   For background on the
implication from the Helson-Szeg\"o condition $PA(1)$ to $(A_2)$, see e.g. Garnett [Gar], Notes
to Ch. VI, Treil and Volberg [TrVo2].\\
\i Under minimality, the relationship between the moving-average coefficients $m = (m_n)$ and
the autoregressive coefficients $r = (r_n)$ becomes symmetrical, and one has the following
complement to Theorem 4:\\

\ni {\bf Theorem 7 (Inoue)}.  For a stationary process $X$, the following are equivalent: \\
(i) The process is minimal. \\
(ii) The autoregressive coefficients $r = (r_n)$ in $(AR)$ satisfy $r \in {\ell}_2$. \\
(iii) $1/h \in H_2$. \\

\ni {\it Proof}.  Since
$$
1/h(z) = \exp \Bigl( \frac{1}{4 \pi} \int \Bigl(\frac{e^{i\t} + z}{e^{i\t} - z}
\Bigr) \log (1/w(\t)) d\t \Bigr) \qquad (z \in D), \eqno(OF')
$$
and $\pm \log w$ are in $L_1$ together, when $1/w \in L_1$ (i.e. the process is minimal) one
can handle $1/w$, $1/h$, $m = (m_n)$ as we handled $w$, $h$ and $r = (r_n)$, giving
$$
1/h \in H_2
$$
and
$$
r = (r_n) \in {\ell}_2.
$$
Conversely, each of these is equivalent to $(min)$; [In1], Prop. 4.2. // \\

\ni {\it 6.4. Rigidity}; $(LM)$, $(CND)$, $(IPF)$. \\
\ni{\it Rigidity; the Levinson-McKean condition}. \\
\i Call $g \in H^1$ {\it rigid} if is determined by its phase or argument:
$$
f \in H^1 \quad \hbox{($f$ not identically 0)}, \quad f/|f| = g/|g| \quad a.e. \quad
\Rightarrow
$$
$$
\quad \hbox{$f = cg$ for some positive constant $c$}.
$$
This terminology is due to Sarason [Sa1], [Sa2]; the alternative terminology, due to Nakazi,
is {\it strongly outer} [Na1], [Na2].  One could instead say that such a function is
{\it determined by its phase}.  The idea originates with de Leeuw and Rudin [dLR] and Levinson
and McKean [LevMcK].  In view of this, we call the condition that $\m$ be absolutely continuous
with spectral density $w = |h|^2$ with $h^2$ rigid, or determined by its phase, the {\it
Levinson-McKean condition}, $(LM)$. \\
\ni {\it Complete non-determinism; intersection of past and future}. \\
\i In [InKa2], the following two conditions are discussed: \\
(i) {\it complete non-determinism},
$$
\hbox{${\cal H}$}_{(-\infty,-1]} \cap \hbox{${\cal H}$}_{[0,\infty)} = \{ 0 \} \eqno(CND)
$$
(for background on this, see [BlJeHa], [JeBl], [JeBlBa]), \\
(ii) the {\it intersection of past and future} property,
$$
\hbox{${\cal H}$}_{(-\infty,-1]} \cap \hbox{${\cal H}$}_{[-n,\infty)} = \hbox{${\cal H}$}_{[-n,-1]} \qquad (n = 1,2,\ldots)
\eqno(IPF)
$$
These are shown to be equivalent in [InKa2].  In [KaBi], it is shown that both are equivalent to
the Levinson-McKean condition, or rigidity:
$$
(LM) \quad \Leftrightarrow \quad (IPF) \quad \Leftrightarrow \quad (CND).
$$
These are weaker than pure minimality ([Bl3], \S 7, [KaBi]).  But since $(CND)$ was already known
to be equivalent to $(PND) \ + \ (IPF)$, they are stronger than $(PND)$.  This takes us from
the weakest of the four intermediate conditions of this section to the stronger of the weak
conditions of \S 3.\\

\ni {\bf 7.  Remarks}\\

\ni 1.  $VMO \subset BMO$.\\
\i The spectral characterizations given above were mainly obtained before the work of Fefferman
[Fe] in 1971, Fefferman and Stein [FeSt] in 1972 (see Garnett [Gar], Ch. VI for a textbook account):
in particular, they predate the {\it Fefferman-Stein decomposition} of a function of bounded
mean oscillation,  $f \in BMO$, as
$$
f = u + \tilde v, \qquad u, v \in L_{\infty}.
$$
This has a complement due to Sarason [Sa3], where $f$ here is in $VMO$ iff $u$, $v$ are
continuous.  Sarason also gives ([Sa3], Th. 2) a characterization of his class of functions of
vanishing mean oscillation $VMO$ within $BMO$ related to Muckenhoupt's condition $(A_2)$. \\
\i While both components $u$, $v$ are needed here, and may be large in norm, it is important to
note that the burden of being large in norm may be born by a {\it continuous} function, leaving
$u$ and $\tilde v$ together to be small in ($L_{\infty}$) norm (in particular, less than
$\pi/2$).  This is the Ibragimov-Rozanov result ([IbRo], V.2 Th. 3), used in \S 6.1 to  show that
absolute regularity (\S 5) implies complete regularity. \\

\ni 2. $H^{1/2} \subset VMO$. \\
\i The class $H^{1/2}$ is contained densely within $VMO$ (Prop. A2, Boutet de Monvel-Berthier
et al. [BouGePu]).  For $H^{1/2}$, one has a version of the Fefferman-Stein decomposition for BMO:
$$
f \in H^{1/2} \quad \Leftrightarrow \quad f = u + \tilde v, \quad u, v \in H^{1/2} \cap
L^{\infty}
$$
([Pel] \S 7.13).\\

\ni 3. {\it Winding number and index}. \\
\i The class $H^{1/2}$ occurs in recent work on {\it topological degree} and
{winding number}; see Brezis [Bre], Bourgain and Kozma [BouKo].  The winding number also occurs
in operator theory as an index in applications of Banach-algebra methods and the Gelfand
transform; see e.g. [Si4], Ch. 5 (cf. Tsirelson [Ts]).\\

\ni 4.  {\it Conformal mapping}. \\
\i The class $H^{1/2}$ also occurs in work of Zygmund on conformal mapping ([Z2], VII.10).\\

\ni 5.  {\it Rapid decay and continuability}. \\
\i Even stronger than the strong conditions considered here in \S \S 4, 5 is assuming that
the Verblunsky coefficients are rapidly decreasing.  This is connected to analytic
continuability of the Szeg\"o function beyond the unit disk; see [Si7]. \\

\ni 6.  {\it Scattering theory}. \\
\i The implication from the strong Szeg\"o (or Golinskii-Ibragimov) condition to the
Helson-Szeg\"o/Helson-Sarason condition $(PA)$ has a recent analogue in scattering theory
(Golinskii et al. [GolKhPeYu], under '$(GI)$ implies $(HS)$'). \\

\ni 7. {\it Wavelets}.\\
\i Traditionally, the subject of time series seemed to consist of two non-intercommunicating
parts, 'time domain' and 'frequency domain' (known to be equivalent to each other via the Kolmogorov Isomorphism Theorem of \S 2).
The subject seemed to suffer from schizophrenia (see e.g. [BriKri] and [HaKR]) -- though the constant relevance of the spectral or frequency side to questions involving time
directly is well illustrated in the apt title 'Past and future' of the paper by Helson and
Sarason [HeSa] (cf. [Pel] \S 8.6).  This unfortunate schism has been healed by the introduction
of {\it wavelet} methods (see e.g. the standard work Meyer [Me], Meyer and Coifman [MeCo], and
in OPUC, Treil and Volberg [TrVo1]).  The practical importance of this may be seen in the
digitization of the FBI's finger-print data-bank (without which the US criminal justice system
would long ago have collapsed).  Dealing with time and frequency together is also crucial in
other areas, e.g. in the high-quality reproduction of classical music. \\

\ni 8. {\it Higher dimensions: matrix OPUC (MOPUC)}. \\
\i We present the theory here in one dimension for simplicity, reserving the case of higher dimensions for a sequel [Bin2].  We note here that
in higher dimensions the measure $\mu$ and the Verblunsky coefficients ${\a}_n$ become matrix-valued (matrix OPUC, or MOPUC), so one loses commutativity.  The multidimensional case
is needed for {\it portfolio theory} in mathematical finance, where one holds a (preferably
balanced) portfolio of risky assets rather than one; see e.g. [BinFrKi]. \\

\ni 9. {\it Non-commutativity}. \\
Much of the theory presented here has a non-commutative analogue in operator theory; see
Blecher and Labuschagne [BlLa], Bekjan and Xu [BeXu] and the references cited there. \\

\ni 10.  {\it Non-stationarity}. \\
\i As mentioned in \S 1, the question of whether or not the process is stationary is vitally important, and
stationarity is a strong assumption.  The basic Kolmogorov Isomorphism Theorem can be extended beyond the stationary case in various ways, e.g. to harmonisable processes (see e.g. [Rao]).  For background, and applications to filtering theory, see e.g. [Kak]; for filtering theory, we refer to e.g. [BaiCr]. \\

\ni 11.  {\it Continuous time}. \\
\i The Szeg\"o condition $(Sz)$ for the unit circle (regarded as the boundary of the unit disc) corresponds to the condition
$$
\int_{-\infty}^{\infty} \frac{\log |f(x)|}{1 + x^2} dx > -\infty
$$
for the real line (regarded as the boundary of the upper half-plane).  This follows from the M\"obius function $w = (z-i)/(z+i)$ mapping the half-plane conformally onto the disc; see e.g. [Du], 189-190.  The consequences of this condition are explored at length in Koosis' monograph on the `logarithmic integral', [Koo2].  Passing from the disc to the half-plane corresponds probabilistically to passing from discrete to continuous time (and analytically to passing from Fourier series to Fourier integrals).  The probabilistic theory is considered at length in Dym and McKean [DymMcK].\\

\ni 12.  {\it Gaussianity and linearity}.\\
\i We have mentioned the close links between Gaussianity
and linearity in \S 1.  For background on Gaussian Hilbert spaces and Fock space, see Janson [Jan],
Peller [Pel]; for extensions to \S \S 5.1, 6 in the Gaussian case, see [IbRo], [Pel], [Bra1] \S 5.  To return to the undergraduate level of our opening paragraph: for an account of
Gaussianity, linearity and regression, see e.g. Williams [Wil], Ch. 8, or [BinFr]. \\

\ni {\bf Acknowledgements}.  This work arises out of collaboration with Akihiko Inoue of Hiroshima University and Yukio Kasahara of Hokkaido University.
It is a pleasure to thank them both.  It is also a pleasure to thank the Mathematics Departments of both universities for their hospitality,  and a Japanese
Government grant for financial support.  I am very grateful to the referee for a thorough and constructive report, which led to many improvements. \\

\begin{center}
{\bf References}
\end{center}
\ni [Ach] N. I. Achieser, {\sl Theory of approximation}, Frederick Ungar, New York, 1956. \\
\ni [BaiCr] A. Bain and D. Crisan, {\sl Fundamentals of stochastic filtering}, Springer, 2009.\\
\ni [BarN-S] O. E. Barndorff-Nielsen and G. Schou, On the parametrization of autoregressive models
by partial autocorrelation.  {\sl J. Multivariate Analysis} {\bf 3} (1973), 408-419.\\
\ni [BasW] E. L. Basor and H. Widom, On a Toeplitz determinant identity of Borodin and Okounkov.
{\sl Integral Equations Operator Theory} {\bf 37} (2000), 397-401.\\
\ni [Bax1] G. Baxter, A convergence equivalence related to polynomials orthogonal on the unit
circle.  {\sl Trans. Amer. Math. Soc.} {\bf 99} (1961), 471-487. \\
\ni [Bax2] G. Baxter, An asymptotic result for the finite predictor.  {\sl Math. Scand.} {\bf 10}
(1962), 137-144. \\
\ni [Bax3] G. Baxter, A norm inequality for a "finite-section" Wiener-Hopf equation.  {\sl
Illinois J. Math.} {\bf 7} (1963), 97-103. \\
\ni [BekXu] T. Bekjan and Q. Xu, Riesz and Szeg\"o type factorizations for non-commutative Hardy
spaces.  {\sl J. Operator Theory} {\bf 62} (2009), 215-231. \\
\ni [Ber] J. Beran, {\sl Statistics for long-memory processes}.  Chapman and Hall, London,
1994.\\
\ni [Berk] K. N. Berk, Consistent autoregressive spectral estimates.  {\sl Ann. Statist.} {\bf 2}
(1974), 489-502. \\
\ni [Beu] A. Beurling, On two problems concerning linear transformations in Hilbert space,  {\sl Acta Math.} {\bf 81} (1948), 239-255
(reprinted in {\sl The collected works of Arne Beurling}, Volumes 1,2, Birkh\"auser, 1989).\\
\ni [Bin1] N. H. Bingham, J\'ozef Marcinkiewicz: Analysis and probability.  {\sl Proc. J\'ozef Marcinkiewicz Centenary Conference} (Pozna\'n, 2010),
Banach Centre Publications {\bf 95} (2011), 27-44.\\
\ni [Bin2] N. H. Bingham: Multivariate prediction and matrix Szeg\"o theory.  Preprint, Imperial College. \\
\ni [BinFr] N. H. Bingham and J. M. Fry, {\sl Regression: Linear models in statistics}.  Springer Undergraduate Mathematics Series, 2010.\\
\ni [BinFrKi] N. H. Bingham, J. M. Fry and R. Kiesel, Multivariate elliptic processes.
{\sl Statistica Neerlandica} {\bf 64} (2010), 352-366. \\
\ni [BinGT] N. H. Bingham, C. M. Goldie and J. L. Teugels, {\sl Regular variation}, 2nd ed.,
Cambridge University Press, 1989 (1st ed. 1987). \\
\ni [BinIK] N. H. Bingham, A. Inoue and Y. Kasahara, An explicit representation of Verblunsky coefficients. {\sl Statistics and Probability Letters},
to appear; online at http://dx.doi.org/10.1016/j.spl.2011.11.004.\\
\ni [BlLa] D. P. Blecher and L. E. Labuschagne, Applications of the Fuglede-Kadison determinant:
Szeg\"o's theorem and outers for non-commutative $H^p$.  {\sl Trans. Amer. Math. Soc.} {\bf 360}
(2008), 6131-6147.\\
\ni [Bl1] P. Bloomfield, An exponential model for the spectrum of a scalar time series.
{\sl Biometrika} {\bf 60} (1973), 217-226. \\
\ni [Bl2] P. Bloomfield, {\sl Fourier analysis of time series: An introduction}.  Wiley, 1976.\\
\ni [Bl3] P. Bloomfield, Non-singularity and asymptotic independence.  {\sl E. J. Hannan
Festschrift, J. Appl. Probab.} {\bf 23A} (1986), 9-21. \\
\ni [BlJeHa] P. Bloomfield, N. P. Jewell and E. Hayashi, Characterization of completely
nondeterministic stochastic processes.  {\sl Pacific J. Math.} {\bf 107} (1983), 307-317.\\
\ni [BogHeTu] B. P. Bogert, M. J. R. Healy and J. W. Tukey, The quefrency alanysis of time series
for echoes: cepstrum, pseudo-autocovariance, cross-cepstrum and saphe cracking.  {\sl Proc.
Symposium on Time Series Analysis} (ed. M. Rosenblatt) Ch. 15, 209-243, Wiley, 1963.\\
\ni [BorOk] A. M. Borodin and A. Okounkov, A Fredholm
determinant formula for Toeplitz
determinants.  {\sl Integral Equations and Operator Theory} {\bf 37} (2000), 386-396. \\
\ni [Bot] A. B\"ottcher, Featured review of the Borodin-Okounkov and Basor-Widom papers.
{\sl Mathematical Reviews} 1790118/6 (2001g:47042a,b). \\
\ni [BotKaSi] A. B\"ottcher, A. Karlovich and B. Silbermann, Generalized Krein algebras and
asymptotics of Toeplitz determinants.  {\sl Methods Funct. Anal. Topology} {\bf 13} (2007), 236-261. \\
\ni [BotSi1] A. B\"ottcher and B. Silbermann, {\sl Analysis of Toeplitz operators}.  Springer,
1990 (2nd ed., with A. Karlovich, 2006). \\
\ni [BotSi2] A. B\"ottcher and B. Silbermann, {\it Introduction to large truncated Toeplitz
matrices}.  Universitext, Springer, 1999. \\
\ni [BotW] A. B\"ottcher and H. Widom, Szeg\"o via Jacobi.  {\sl Linear Algebra and Applications}
{\bf 419} (2006), 656-667.\\
\ni [BouKo] J. Bourgain and G. Kozma, One cannot hear the winding number. {\sl J. European Math.
Soc.} {\bf 9} (2007), 637-658.\\
\ni [BoGePu] A. Boutet de Monvel-Berthier, V. Georgescu and R. Purice, A boundary-value problem
related to the Ginzburg-Landau model.  {\sl Comm. Math. Phys.} {\bf 142} (1991), 1-23. \\
\ni [BoxJeRe] G. E. P. Box, G. M. Jenkins and G. C. Reinsel, {\sl Time-series analysis.
Forecasting and control} (4th ed.).  Wiley, 2008.\\
\ni [Bra1] R. C. Bradley, Basic properties of strong mixing conditions.  Pages 165-192 in [EbTa].\\
\ni [Bra2] R. C. Bradley, Basic properties of strong mixing conditions.  A survey and some open
questions.  {\sl Probability Surveys} {\bf 2} (2005), 107-144. \\
\ni [Bra3] R. C. Bradley, {\sl Introduction to strong mixing conditions}, Volumes 1-3.  Kendrick
Press, Heber City, UT, 2007. \\
\ni [dBR] L. de Branges and J. Rovnyak, {\sl Square-summable power series}.  Holt, Rinehart and
Winston, New York, 1966. \\
\ni [Bre] H. Brezis, New questions related to the topological
degree. {\sl The unity of mathematics}. 137-154, Progr. Math. {\bf 244} (2006), Birkh\"auser,
Boston MA. \\
\ni [Bri] D. R. Brillinger, John W. Tukey's work on time series
and spectrum analysis.  {\sl Ann. Statist.} {\bf 30} (2002), 1595-1918. \\
\ni [BriKri] D. R. Brillinger and P. R. Krishnaiah (ed.), {\sl Time series in the frequency domain}.  Handbook of Statistics {\bf 3}, North-Holland, 1983. \\
\ni [BroDav] P. J. Brockwell and R. A. Davis, {\sl Time series: Theory and methods}
(2nd ed.), Springer, New York, 1991 (1st ed. 1987).\\
\ni [Cra] H. Cram\'er, On harmonic analysis in certain function spaces.  {\sl Ark. Mat. Astr. Fys.} {\bf 28B} (1942), 1-7,
reprinted in {\sl Collected Works of Harald Cram\'er} Volume II, 941-947, Springer, 1994. \\
\ni [CraLea] H. Cram\'er and R. Leadbetter, {\sl Stationary and related stochastic processes}.  Wiley, 1967. \\
\ni [Ch] K.-L. Chung, {\sl A course in probability theory}, 3rd ed.  Academic Press, 2001
(2nd ed. 1974, 1st ed. 1968). \\
\ni [CorFoSi] I. P. Cornfeld, S. V. Fomin and Ya. G. Sinai, {\sl Ergodic theory}.  Grundl. math.
Wiss. {\bf 245}, Springer, 1982. \\
\ni [CovTh] T. M. Cover and J. A. Thomas, {\sl Elements of information theory}.  Wiley, 1991. \\
\ni [Cox] D. R. Cox, Long-range dependence: a review.  Pages 55-74 in {\sl Statistics: An
appraisal} (ed. H. A. David and H. T. David), Iowa State University Press, Ames IA, reprinted in
{\sl Selected statistical papers of Sir David Cox, Volume 2}, 379-398, Cambridge University
Press, 2005. \\
\ni [Deb] L. Debowski, On processes with summable partial autocorrelations.  {\sl Statistics and Probability Letters} {\bf 77} (2007), 752-759. \\
\ni [Deg] S. D\'egerine, Canonical partial autocorrelation function of a multivariate time
series.  {\sl Ann. Statist.} {\bf 18} (1990), 961-971.\\
\ni [dLR] K. de Leeuw and W. Rudin, Extreme points and extremum problems in $H^1$.  {\sl Pacific
J. Math.} {\bf 8} (1958), 467-485.\\
\ni [Dia] P. Diaconis, G. H. Hardy and probability??? {\sl Bull. London Math. Soc.} {\bf 34}
(2002), 385-402.\\
\ni [Do] J. L. Doob, {\sl Stochastic processes}.  Wiley, 1953. \\
\ni [DouOpTa] P. Doukhan, G. Oppenheim and M. S. Taqqu (ed.), {\sl Theory and applications of long-
range dependence}.  Birkh\"auser, Basel, 2003. \\
\ni [DunSch] N. Dunford and J. T. Schwartz, {\sl Linear operators, Part II:
Spectral theory: Self-adjoint operators on Hilbert space}.  Interscience, 1963. \\
\ni [Dur] J. Durbin, The fitting of time-series models.  {\sl Rev. Int. Statist. Inst.} {\bf 28}
(1960), 233-244. \\
\ni [Du] P. L. Duren, {\sl Theory of} $H^p$ {\sl spaces}.  Academic Press, New York, 1970. \\
\ni [Dym] H. Dym, M. G. Krein's contributions to prediction theory.  {\sl Operator Theory Advances and Applications} {\bf 118} (2000), 1-15, Birkh\"auser, Basel. \\
\ni [DymMcK] H. Dym and H. P. McKean, {\sl Gaussian processes, function theory and the inverse spectral problem}.  Academic Press, 1976.\\
\ni [EbTa] E. Eberlein and M. S. Taqqu (ed.), {\sl Dependence in probability and statistics.  A
survey of recent results}.  Birkh\"auser, 1986. \\
\ni [Ell] R. S. Ellis, {\sl Entropy, large deviations and statistical mechanics}.  Grundl. math. Wiss {\bf 271}, Springer, 1985.\\
\ni [Fe] C. Fefferman, Characterizations of bounded mean oscillation.  {\sl Bull. Amer. Math. Soc.} {\bf 77} (1971), 587-588.\\
\ni [FeSt] C. Fefferman and E. M. Stein, $H^p$ spaces of several variables.  {\sl Acta Math.} {\bf 129} (1972), 137-193.\\
\ni [Gao] J. Gao, {\sl Non-linear time series.  Semi-parametric and non-parametric models}.
Monogr. Stat. Appl. Prob. {\bf 108}, Chapman and Hall, 2007.\\
\ni [Gar] J. B. Garnett, {\sl Bounded analytic functions}.  Academic Press, 1981 (Grad. Texts in Math. {\bf 236}, Springer, 2007).\\
\ni [Geo] H.-O. Georgii, {\sl Gibbs measures and phase transitions}.  Walter de Gruyter, 1988. \\
\ni [GerCa] J. S. Geronimo and K. M. Case, Scattering theory and polynomials orthogonal on the
unit circle.  {\sl J. Math. Phys.} {\bf 20} (1979), 299-310. \\
\ni [Ges] F. Gesztesy, P. Deift, C. Galves, P. Perry and W. Schlag (ed.), {\sl Spectral theory and mathematical physics:
A Festschrift in honor of Barry Simon's sixtieth birthday}.  Proc. Symp. Pure Math. {\bf 76} Parts 1, 2, Amer. Math. Soc., 2007. \\
\ni [GiKoSu] L. Giraitis, H. L. Koul and D. Surgailis, {\sl Large sample inference for long memory processes}, World Scientific, 2011. \\
\ni [GolKPY] L. Golinskii, A. Kheifets, F. Peherstorfer and P. Yuditskii, Scattering theory for CMV matrices: uniqueness, Helson-Szeg\"o and strong Szeg\"o theorems.  {\sl Integral Equations and Operator Theory} {\bf 69} (2011), 479-508.\\
\ni [GolTo] L. Golinskii and V. Totik, Orthogonal polynomials: from Jacobi to Simon.  P. 715-742 in [Ge], Part 2. \\
\ni [GolvL] G. H. Golub and C. F. van Loan, {\sl Matrix computations}, 3rd ed.,
Johns Hopkins University Press, 1996 (1st ed. 1983, 2nd ed. 1989). \\
\ni [GrSz] U. Grenander and G. Szeg\"o, {\sl Toeplitz forms and their applications}.  University
of California Press, Berkeley CA, 1958. \\
\ni [Gri1] G. R. Grimmett, {\sl Percolation}, 2nd ed.  Grundl. math. Wiss. {\bf 321}, Springer,
1999 (1st ed. 1989).\\
\ni [Gri2] G. R. Grimmett, {\sl The random cluster model}.  Grundl. math. Wiss. {\bf 333},
Springer, 2006.\\
\ni [Ha1] E. J. Hannan, {\sl Multiple time series}.  Wiley, 1970. \\
\ni [Ha2] E. J. Hannan, The Whittle likelihood and frequency estimation.  Chapter 15 (p. 205-212) in [Kel]. \\
\ni [HaKR] E. J. Hannan, P. K. Krishnaiah and M. M. Rao, {\sl Time series in the time domain}.  Handbook of Statistics {\bf 5}, North-Holland, 1985. \\
\ni [He] H. Helson, {\sl Harmonic analysis}, 2nd ed., Hindustan Book Agency, 1995.\\
\ni [HeSa] H. Helson and D. Sarason, Past and future.  {\sl Math. Scand} {\bf 21} (1967),
5-16.\\
\ni [HeSz] H. Helson and G. Szeg\"o, A problem in prediction theory.  {\sl Acta Mat. Pura Appl.}
{\bf 51} (1960), 107-138.\\
\ni [Ho] K. Hoffman, {\sl Banach spaces of analytic functions}, Prentice-Hall, Englewood Cliffs
NJ, 1962. \\
\ni [Hos] J. R. Hosking, Fractional differencing.  {\sl Biometrika} {\bf 68} (1981), 165-176. \\
\ni [HuMuWe] R. A. Hunt, B. Muckenhoupt and R. L. Wheeden, Weighted norm inequalities for the conjugate
function and Hilbert transform.  {\sl Trans. Amer. Math. Soc.} {\bf 176} (1973), 227-151. \\
\ni [IbLi] I. A. Ibragimov and Yu. V. Linnik, {\sl Independent and stationary sequences of random
variables}.  Wolters-Noordhoff, 1971. \\
\ni [IbRo] I. A. Ibragimov and Yu. A. Rozanov, {\sl Gaussian random processes}.  Springer, 1978.\\
\ni [In1] A. Inoue, Asymptotics for the partial autocorrelation function of a stationary process.
{\sl J. Analyse Math.} {\bf 81} (2000), 65-109. \\
\ni [In2] A. Inoue, Asymptotic behaviour for partial autocorrelation functions of fractional
ARIMA processes.  {\sl Ann. Appl. Probab.} {\bf 12} (2002), 1471-1491. \\
\ni [In3] A. Inoue, AR and MA representations of partial autocorrelation functions, with
applications.  {\sl Prob. Th. Rel. Fields} {140} (2008), 523-551. \\
\ni [InKa1] A. Inoue and Y. Kasahara, Partial autocorrelation functions of the fractional ARIMA
processes.  {\sl J. Multivariate Analysis} {\bf 89} (2004), 135-147. \\
\ni [InKa2] A. Inoue and Y. Kasahara, Explicit representation of finite predictor coefficients
and its applications.  {\sl Ann. Statist.} {\bf 34} (2006), 973-993. \\
\ni [Jan] S. Janson, {\sl Gaussian Hilbert spaces}.  Cambridge Tracts in Math. {\bf 129}, Cambridge
University Press, 1997.\\
\ni [JeBl] N. P. Jewell and P. Bloomfield, Canonical correlations of past and future for time
series: definitions and theory.  {\sl Ann. Statist.} {\bf 11} (1983), 837-847.\\
\ni [JeBlBa] N. P. Jewell, P. Bloomfield and F. C. Bartmann, Canonical correlations of past and
future for time series: bounds and computation.  {\sl Ann. Statist.} {\bf 11} (1983), 848-855.\\
\ni [Kac] M. Kac, Toeplitz matrices, transition kernels and a related problem in probability
theory.  {\sl Duke Math. J.} {\bf 21} (1954), 501-509. \\
\ni [KahSa] J.-P. Kahane and R. Salem, {\sl Ensembles parfaits et s\'eries trigonom\'etriques}, 2nd ed.
Hermann, Paris, 1994. \\
\ni [Kak] Y. Kakihara, The Kolmogorov isomorphism theorem and extensions to some non-stationary processes.
{\sl Stochastic processes: Theory and methods} (ed. D. N. Shanbhag and C. R. Rao),
{\sl Handbook of Statistics} {\bf 19}, North-Holland, 2001, 443-470.\\
\ni [KanSch] H. Kantz and T. Schreiber, {\sl Nonlinear time series analysis}, Cambridge University Press, 1997 (2nd ed. 2004). \\
\ni [KaBi] Y. Kasahara and N. H. Bingham, Verblunsky coefficients and Nehari sequences.  Preprint, Hokkaido University.\\
\ni [KatSeTe] D. Kateb, A. Seghier and G. Teyssi\`ere, Prediction, orthogonal polynomials and
Toeplitz matrices.  A fast and reliable approach to the Durbin-Levinson algorithm.  Pages
239-261 in [TK].\\
\ni [Kel] F. P. Kelly (ed.), {\sl Probability, statistics and optimization.  A tribute to Peter
Whittle}.  Wiley, 1994. \\
\ni [KenSt] M. G. Kendall and A. Stuart, {\sl The advanced theory of statistics}.  Charles Griffin.
Volume 1 (4th ed., 1977), Vol. 2 (3rd ed, 1973), vol. 3 (3rd ed., 1976).\\
\ni [KokTa] P. S. Kokoszka and M. S. Taqqu, Can one use the Durbin-Levinson algorithm to generate
infinite-variance fractional ARIMA time series?  {\sl J. Time Series Analysis} {\bf 22} (2001),
317-337. \\
\ni [Kol] A. N. Kolmogorov, Stationary sequences in Hilbert space.  {\sl Bull. Moskov. Gos. Univ.
Mat.} {\bf 2} (1941), 1-40 (in Russian; reprinted, {\sl Selected works of A. N. Kolmogorov,
Vol. 2: Theory of probability and mathematical statistics}, Nauka, Moskva, 1986, 215-255). \\
\ni [Koo1] P. Koosis, {\sl Introduction to} $H^p$ {\sl spaces}, 2nd ed.  Cambridge Tracts Math.
{\bf 115}, Cambridge Univ. Press, 1998 (1st ed. 1980).\\
\ni [Koo2] P. Koosis, {\sl The logarithmic integral}, I, 2nd ed., Cambridge Univ. Press, 1998 (1st ed. 1988), II, Cambridge Univ. Press, 1992.\\
\ni [Kr] M. G. Krein, On some new Banach algebras and Wiener-L\'evy type theorems for Fourier
series and integrals.  {\sl Amer. Math. Soc. Translations} (2) {\bf 93} (1970), 177-199
(Russian original: {\sl Mat. Issled.} {\bf 1} (1966), 163-288).\\
\ni [Lev] N. Levinson, The Wiener (RMS) error criterion in filter design and prediction.
{\sl J. Math. Phys. MIT} {\bf 25} (1947), 261-278. \\
\ni [LevMcK] N. Levinson and H. P. McKean, Weighted trigonometrical approximation on $R^1$ with
application to the germ field of a stationary Gaussian noise.  {\sl Acta Math.} {\bf 112}
(1964), 99-143.\\
\ni [Li] L. M. Li, Some notes on mutual information between past and future.  {\sl J. Time
Series Analysis} {\bf 27} (2006), 309-322. \\
\ni [LiXi] L. M. Li and Z. Xie, Model selection and order determination for time series by
information between the past and the future.  {\sl J. Time Series Analysis} {\bf 17} (1996),
65-84. \\
\ni [LuZhKi] R. Lund, Y. Zhao and P. C. Kiessler, Shapes of stationary autocovariances.  {\sl J.
Applied Probability} {\bf 43} (2006), 1186-1193. \\
\ni [Ly1] R. Lyons, Characterizations of measures whose Fourier-Stieltjes transforms vanish at
infinity.  {\sl Bull. Amer. Math. Soc.} {\bf 10} (1984), 93-96.\\
\ni [Ly2] R. Lyons, Fourier-Stieltjes coefficients and asymptotic distribution modulo 1.
{\sl Ann. Math.} {\bf 122} (1985), 155-170. \\
\ni [Ly3] R. Lyons, Seventy years of Rajchman measures.  {\sl J. Fourier Anal. Appl.},
Kahane Special Issue (1995), 363-377.\\
\ni [MakWe] A. Makagon and A. Weron, $q$-variate minimal stationary processes.  {\sl Studia
Math.} {\bf 59} (1976), 41-52. \\
\ni [MatNeTo] A. M\'at\'e, P. Nevai and V. Totik, Aymptotics for the ratio of leading coefficients
of orthogonal polynomials on the unit circle.  {\sl Constructive Approximation} {\bf 1} (1985),
63-69. \\
\ni [McCW] B. M. McCoy and T. T. Wu, {\sl The two-dimensional Ising model}.  Harvard Univ. Press,
Cambridge MA, 1973. \\
\ni [McC] P. McCullagh, John Wilder Tukey, 1915-2000.  {\sl Biographical Memoirs of Fellows of
the Royal Society} {\bf 49} (2003), 537-555.\\
\ni [McLZ] A. I. McLeod and Y. Zhang, Partial autocorrelation parametrization for subset
regression.  {\sl J. Time Series Analysis} {\bf 27} (2006), 599-612. \\
\ni [Me] Y. Meyer, {\sl Wavelets and operators}.  Cambridge Univ. Press, 1992. \\
\ni [MeCo] Y. Meyer and R. Coifman, {\sl Wavelets.  Calder\'on-Zygmund and multilinear
operators}.  Cambridge Univ. Press, 1997.\\
\ni [Nak1] T. Nakazi, Exposed points and extremal problems in $H^1$.  {\sl J. Functional
Analysis} {\bf 53} (1983), 224-230.\\
\ni [Nak2] T. Nakazi, Exposed points and extremal problems in $H^1$, II.  {\sl Tohoku Math. J.}
{\bf 37} (1985), 265-269.\\
\ni [vN] J. von Neumann, Allgemeine Eigenwerttheorie Hermitescher Funktionaloperatoren.  {\sl
Math. Ann.} {\bf 102}, 49-131 ({\sl Collected Works} II.1). \\
\ni [Nik1] N. K. Nikolskii, {\sl Treatise on the shift operator: Spectral function theory}.
Grundl. math. Wiss. {\bf 273}, Springer, 1986.\\
\ni [Nik2] N. K. Nikolskii, {\sl Operators, functions and systems: an easy reading.  Volume 1:
Hardy, Hankel and Toeplitz; Volume 2: Model operators and systems}.  Math. Surveys and Monographs
{\bf 92, 93}, Amer. Math. Soc., 2002. \\
\ni [OpSc] A. V. Oppenheim and R. W. Schafer, {\sl Discrete signal processing}.  Prentice-Hall,
1989. \\
\ni [Pel] V. V. Peller, {\sl Hankel operators and their applications}.  Springer, 2003. \\
\ni [PoSz] G. P\'olya and G. Szeg\"o, {\sl Problems and theorems in analysis}, I, II.  Classics
in Math., Springer, 1998 (transl. 4th German ed., 1970; 1st ed. 1925).\\
\ni [Pou] M. Pourahmadi, {\sl Foundations of time series analysis and prediction theory}.
Wiley, 2001. \\
\ni [Rak] E. A. Rakhmanov, On the asymptotics of the ratios of orthogonal polynomials, II,
{\sl Math. USSR Sb.} {\bf 58} (1983), 105-117. \\
\ni [Ram] F. L. Ramsey, Characterization of the partial autocorrelation function.  {\sl Ann.
Statist.} {\bf 2} (1974), 1296-1301.\\
\ni [Rao] M. M. Rao, Harmonizable, Cram\'er and Karhunen classes of processes.  Ch. 10 (p.276-310) in [HaKR].\\
\ni [Rob] P. M. Robinson (ed.), {\sl Time series with long memory}.  Advanced Texts in
Econometrics, Oxford University Press, 2003. \\
\ni [RoRo] M. Rosenblum and J. Rovnyak, {\sl Hardy classes and operator theory}, Dover, New York,
1997 (1st ed. Oxford University Press, 1985). \\
\ni [Roz] Yu. A. Rozanov, {\sl Stationary random processes}.  Holden-Day, 1967. \\
\ni [Ru] W. Rudin, {\sl Real and complex analysis}, 2nd ed.  McGraw-Hill, 1974 (1st ed. 1966). \\
\ni [Sa1] D. Sarason, {\sl Function theory on the unit circle}.  Virginia Polytechnic Institute
and State University, Blacksburg VA, 1979. \\
\ni [Sa2] D. Sarason, An addendum to "Past and future", {\sl Math. Scand.} {\bf 30} (1972),
62-64.\\
\ni [Sa3] D. Sarason, Functions of vanishing mean oscillation, {\sl Trans. Amer. Math. Soc.}
{\bf 207} (1975), 391-405.\\
\ni [Si1] B. Simon, {\sl The statistical mechanics of lattice gases, Volume 1}.  Princeton
University Press, 1993.\\
\ni [Si2] B. Simon, The Golinskii-Ibragimov method and a theorem of Damanik-Killip.
{\sl Int. Math. Res. Notes} (2003), 1973-1986. \\
\ni [Si3] B. Simon, OPUC on one foot.  {\sl Bull. Amer. Math. Soc.} {\bf 42} (2005), 431-460.\\
\ni [Si4] B. Simon, {\sl Orthogonal polynomials on the unit circle.  Part 1: Classical theory}.
AMS Colloquium Publications 54.1, American Math. Soc., Providence RI, 2005. \\
\ni [Si5] B. Simon, {\sl Orthogonal polynomials on the unit circle.  Part 2: Spectral theory}.
AMS Colloquium Publications 54.2, American Math. Soc., Providence RI, 2005. \\
\ni [Si6] B. Simon, The sharp form of the strong Szeg\"o theorem.  {\sl Contemorary Math.}
{\bf 387} (2005), 253-275, AMS, Providence RI. \\
\ni [Si7] B. Simon, Meromorphic Szeg\"o functions and asymptotic series for Verblunsky
coefficients.  {\sl Acta Math.} {\bf 195} (2005), 267-285. \\
\ni [Si8] B. Simon, Ed Nelson's work in quantum theory.  {\sl Diffusion, quantum theory and
radically elementary mathematics} (ed. W. G. Faris), {\sl Math. Notes} {\bf 47} (2006), 75-93.\\
\ni [Si9] B. Simon, {\sl Szeg\"o's theorem and its descendants: Spectral theory for} $L^2$ {\sl perturbations of orthogonal polynomials}.  Princeton University Press, 2011. \\
\ni [Si10] B. Simon, {\sl Convexity: Ana analytic viewpoint}.  Cambridge Tracts in Math. {\bf 187}, Cambridge University Press, 2011. \\
\ni [Sz1] G. Szeg\"o, Ein Grenzwertsatz \"uber die Toeplitzschen Determinanten einer reellen
positiven Funktion.  {\sl Math. Ann.} {\bf 76} (1915), 490-503.\\
\ni [Sz2] G. Szeg\"o, Beitr\"age zur Theorie der Toeplitzschen
Formen.  {\sl Math. Z.} {\bf 6} (1920), 167-202. \\
\ni [Sz3] G. Szeg\"o, Beitr\"age zur Theorie der Toeplitzschen
Formen, II.  {\sl Math. Z.} {\bf 9} (1921), 167-190. \\
\ni [Sz4] G. Szeg\"o, {\sl Orthogonal polynomials}.  AMS Colloquium Publications 23, American
Math. Soc., Providence RI, 1939.\\
\ni [Sz5] G. Szeg\"o, On certain Hermitian forms associated with the Fourier series of a
positive function.  {\sl Festschrift Marcel Riesz} 222-238, Lund, 1952. \\
\ni [SzNF] B. Sz.-Nagy and C. Foias, {\sl Harmonic analysis of operators on Hilbert space}, North-Holland, 1970
(2nd ed., with H. Bercovici and L. K\'erchy, Springer Universitext, 2010).\\
\ni [TeKi] G. Teyssi\`ere and A. P. Kirman (ed.), {\sl Long memory in economics}.  Springer, 2007.\\
\ni [Tor] A. Torchinsky, {\sl Real-variable methods in harmonic analysis}, Dover, 2004 (Academic Press, 1981).\\
\ni [TrVo1] S. Treil and A. Volberg, Wavelets and the angle between past and future.  {\sl J. Functional analysis} {\bf 143} (1997), 269-308.\\
\ni [TrVo2] S. Treil and A. Volberg, A simple proof of the Hunt-Muckenhoupt-Wheeden theorem.  Preprint, 1997.\\
\ni [Ts] B. Tsirelson, Spectral densities describing off-white noises.  {\sl Ann. Inst. H. Poincar\'e Prob. Stat.} {\bf 38} (2002), 1059-1069.\\
\ni [V1] On positive harmonic functions.  A contribution to the
algebra of Fourier series.  {\sl Proc. London Math. Soc.} {\bf 38}
(1935), 125-157. \\
\ni [V2] On positive harmonic functions (second paper). {\sl Proc.
London Math. Soc.} {\bf 40} (1936), 290-320. \\
\ni [Wh] P. Whittle, {\sl Hypothesis testing in time series analysis}. Almqvist and
Wiksell, Uppsala, 1951. \\
\ni [Wi1] N. Wiener, Generalized harmonic analysis.  {\sl Acta Math.} {\bf 55} (1930), 117-258
(reprinted in {\sl Generalized harmonic analysis and Tauberian theorems},  MIT Press, Cambridge
MA, 1986, and {\sl Collected Works, Volume II: Generalized harmonic analysis and Tauberian theory; classical harmonic and complex analysis} (ed. P. Masani), MIT Press, Cambridge MA, 1979).\\
\ni [Wi2] N. Wiener, {\sl Extrapolation, interpolation and smoothing of stationary time
series.  With engineering applications}.  MIT Press/Wiley, 1949.\\
\ni [Wi3] N. Wiener, {\sl Collected Works, Volume III: The Hopf-Wiener integral equation; prediction and filtering; quantum mechanics and relativity; miscellaneous mathematical papers} (ed. P. Masani), MIT Press, Cambridge MA, 1981. \\
\ni [Wil] D. Williams, {\sl Weighing the odds}.  Cambridge University Press, 2001.\\
\ni [Wo] H. Wold, {\sl A study in the analysis of stationary time series}.  Almqvist and
Wiksell, Uppsala, 1938 (2nd ed., appendix by Peter Whittle, 1954). \\
\ni [Z1] A. Zygmund, Sur les fonctions conjugu\'ees.  {\sl Fund. Math.} {\bf 13} (1929), 284-303;
corr. {\sl Fund. Math.} {\bf 18} (1932), 312 (reprinted in [Z3], vol 1). \\
\ni [Z2] A. Zygmund, {\sl Trigonometric series}, Volumes 1,2, Cambridge University Press, 1968.\\
\ni [Z3] A. Zygmund, {\sl Selected papers of Antoni Zygmund} (ed. A. Hulanicki, P. Wojtaszczyk and W. Zelasko), Volumes 1-3, Kluwer, Dordrecht, 1989.\\

\ni N. H. Bingham, Mathematics Department, Imperial College London, London SW7 2AZ, UK $\quad$
nick.bingham@btinternet.com $\quad$ n.bingham@ic.ac.uk \\

\end{document}